\documentclass[10pt]{article}
 \usepackage{benstyle}
 
\DefineNamedColor{named}{orange} {rgb}{1,0.55,0}

\title{Beyond consistent reconstructions: optimality and sharp bounds for generalized sampling, and application to the uniform resampling problem}
\author{Ben Adcock \\ Department of Mathematics \\ Purdue University \\ \hspace{20pc} \\ \vspace{-2.25pc} \\
150 N. University Street \\ West Lafayette, IN 47907 \\ USA  \and Anders C. Hansen \\ DAMTP, Centre for Mathematical Sciences \\ University of Cambridge \\ Wilberforce Rd, Cambridge CB3 0WA \\ United Kingdom \and Clarice Poon \\ DAMTP, Centre for Mathematical Sciences \\ University of Cambridge \\ Wilberforce Rd, Cambridge CB3 0WA \\ United Kingdom }
\begin{document}
\maketitle

\begin{abstract}
Generalized sampling is a recently developed linear framework for sampling and reconstruction in separable Hilbert spaces.  It allows one to recover any element in any finite-dimensional subspace given finitely many of its samples with respect to an arbitrary frame.  Unlike more common approaches for this problem, such as the consistent reconstruction technique of Eldar et al, it leads to completely stable numerical methods possessing both guaranteed stability and accuracy.

The purpose of this paper is twofold.  First, we give a complete and formal analysis of generalized sampling, the main result of which being the derivation of new, sharp bounds for the accuracy and stability of this approach.  Such bounds improve those given previously, and result in a necessary and sufficient condition, the stable sampling rate, which guarantees \textit{a priori} a good reconstruction.  Second, we address the topic of optimality.  Under some assumptions, we show that generalized sampling is an optimal, stable reconstruction.  Correspondingly, whenever these assumptions hold, the stable sampling rate is a universal quantity.  In the final part of the paper we illustrate our results by applying generalized sampling to the so-called uniform resampling problem.
\end{abstract}

\section{Introduction}\label{s:introduction}
A central theme in sampling theory is the recovery of a signal or an image from a collection of its measurements.  Mathematically, this can be modelled in a separable Hilbert space $\rH$, with the samples of the unknown signal $f \in \rH$ being of the form
\bes{
\hat{f}_{j} = \ip{f}{\psi_{j}},\quad j=1,2,\ldots ,
}
where $\{ \psi_j \}^{\infty}_{j=1}$ is a collection of elements belonging to $\rH$ (here $\ip{\cdot}{\cdot}$ is the inner product on $\rH$).  Typically, the \textit{sampling} system $\{ \psi_j \}^{\infty}_{j=1}$ forms a frame for its span $\rS = \overline{\spn \{ \psi_1,\psi_2,\ldots \}}$.

One of the most common, and arguably one of the most important, examples of this type of sampling is the recovery of a function $f$ with compact support from pointwise evaluations of its Fourier transform $\hat{f}$.  In this case, $\rH = \rL^2(-1,1)^d$, where $\mathrm{supp}(f) \subseteq (-1,1)^d$ without loss of generality, and $\psi_j(x)=\E^{\I \pi \omega_j \cdot x}$ for suitable values $\{ \omega_j \}_{j \in \bbN} \subseteq \bbR^d$.  This is precisely the type of sampling encountered in Magnetic Resonance Imaging (MRI), for example.  

If the measurements $\{ \omega_j \}_{j \in \bbN}$ are taken uniformly, both $f$ and $\hat{f}$ can be recovered via the Shannon Sampling Theorem \cite{jerrishannon,unser2000sampling}.    However, the slow convergence of the corresponding reconstructions (both infinite sums), as well as the appearance of the Gibbs phenomenon, means that this approach is often not practical \cite{EldarMag,ParkerKenyonImage,unser2000sampling}.  In cases where measurements $\{ \omega_j \}_{j \in \bbN}$ are not uniformly distributed, no simple reconstruction need exist.  It is standard in this setting to use a gridding algorithm \cite{JacksonEtAlGridding,MeyerEtAllGridding,RosenfeldURS1,GelbNonuniformFourier}.  However, this also typically leads to less than satisfactory accuracy.  Much as in the uniform case, unsightly Gibbs oscillations also persist \cite{GelbNonuniformFourier}.

The MRI problem serves to illustrate several key issues that are critical to this paper.  First, although $f$ is sampled via an infinite collection of elements $\{ \psi_j \}^{\infty}_{j=1}$, in practice we only have access to a finite number.  Thus, the problem we consider throughout this paper is that of recovering $f$ from only its \textit{first $n$} samples $\hat{f}_1,\ldots,\hat{f}_n$.  Key issues herein are those of \textit{approximation} -- namely, how well $f$ can be recovered as $n \rightarrow \infty$ -- and \textit{robustness} -- does increasing $n$ lead to worse stability, and thus more sensitivity to noise and round-off error?  

The MRI example also highlights another key issue.  Namely, the samples $\{\hat{f}_j \}_{j \in \bbN}$ of $f$ are fixed, and cannot easily be altered. This situation occurs typically when the sampling scheme is specified by some physical device, e.g.\ the MR scanner in the above example.  Although it is actually possible to modify MR scanners to acquire different types of measurements, such as \textit{wavelet-encoded} MRI \cite{HealyWeaverWaveletIEEE,WeaverEtAlWaveletEncoding}, this is not without complications \cite{LaineWaveletsBiomed}.  Thus, the question we consider in this paper is the following: given a finite number of fixed samples of an element $f$ of a Hilbert space $\rH$, how can one obtain a good (i.e.\ accurate and robust) reconstruction?

This question is not new, and there has been much interest in the last several decades in alternative reconstructions to those given by the Shannon Sampling Theorem.  This is typically based on the following principle: many signals that arise in practice can be much better represented in terms of a different collection of elements $\{ \phi_j \}_{j \in \bbN} \subseteq \rH$ \cite{EldarMag,unser2000sampling} than by Shannon's theorem.  Common examples of such systems include wavelets, splines and polynomials, as well as more exotic objects such as curvelets \cite{Cand,candes2004new}, shearlets \cite{Gitta,Gitta2,Gitta3} and contourlets \cite{Vetterli,Do}.  Thus, given this additional knowledge about $f$, the problem is now as follows: how can we compute a reconstruction in the system $\{ \phi_j \}_{j \in \bbN}$ from the measurements $\hat{f}_{j} = \ip{f}{\psi_j}$?  

Consistent sampling is a linear technique designed specifically for this problem, based on stipulating that the reconstruction of $f$ agrees with the available measurements.  Introduced by Unser \& Aldroubi \cite{unser1994general,unserzerubia} and later generalised significantly by Eldar et al \cite{EldarRobConsistSamp,eldar2003FAA,eldar2003sampling,eldar2005general}, this technique has proved successful in a number of areas, and is quite widely used in practice \cite{unser2000sampling}.  However, there are a number of drawbacks.  As discussed in \cite{BAACHShannon,BAACHSampTA,EldarMinimax,UnserHirabayashiConsist}, consistent reconstructions need not be stable or convergent as the number of measurements increases.  Whilst stability and convergence can be guaranteed in certain shift-invariant spaces  \cite{eldar2003sampling,UnserZerubiaGS}, it is quite easy to devise examples outside of this setting for which consistent reconstructions either fail to converge, or are extremely unstable, or both \cite{BAACHShannon}.  

Fortunately, it transpires that these issues can be overcome by using an recently-introduced alternative technique, known as \textit{generalized sampling} \cite{BAACHShannon,BAACHSampTA,BAACHAccRecov}.  This method forms the primary focus of our paper.  Our main results are described in the next section.

\subsection{Novelty of the paper and overview}
The purpose of this paper is to give a complete analysis of generalized sampling.  Not only do we establish sharp bounds for stability and the reconstruction error, which improve on those appearing previously in \cite{BAACHShannon,BAACHSampTA,BAACHAccRecov}, we also provide several optimality results.  These results demonstrate that under mild conditions generalized sampling cannot be outperformed.  To illustrate our results we consider the so-called \textit{uniform resampling} problem.

Let us now give an overview of the remainder of the paper.  In \S \ref{s:reconprob}--\ref{s:consistent} we give a mathematical description of the general instability and nonconvergence of consistent reconstructions.  For this, it is necessary to present a formal description of the reconstruction problem (\S \ref{s:reconprob}).  Herein we also introduce two key constants to assess different methods for the problem: the \textit{condition number} $\kappa$ and the \textit{quasi-optimality constant} $\mu$.  The former measures the sensitivity of a given reconstruction method to perturbations (e.g. noise and round-off), whilst the latter quantifies how close the reconstruction of an element $f$ is to its best (i.e.\ optimal) reconstruction in the desired system of functions.  For succinctness, we also introduce the reconstruction constant $C$ of a method, defined as the maximum of $\kappa$ and $\mu$.  In \S \ref{s:reconprob} we explain why it is vital in practice to have a small reconstruction constant $C$.

The focus of \S \ref{s:consistent} is consistent sampling.  By analysing the reconstruction constant $C$ in this instance, we provide a comprehensive answer as to when this approach will give poor (i.e.\ unstable and inaccurate) reconstructions.  Moreover, we show how one can determine \textit{a priori} an answer to this question by performing a straightforward computation.  In other words, the success or failure of a consistent reconstruction for a particular problem can always be determined beforehand.  We also show how this question can be reinterpreted in terms of the behaviour \textit{finite sections} of infinite operators, and thus draw a connection between problems in sampling and computational spectral theory (such a connection was first discussed in \cite{BAACHShannon}).

To overcome the issues inherent to consistent reconstructions the new approach of generalized sampling was introduced in \cite{BAACHShannon,BAACHSampTA,BAACHAccRecov}.  In the second part of the paper (\S \ref{s:GS}--\ref{ss:stabsamp}) we improve the previous analysis of \cite{BAACHShannon,BAACHAccRecov} by using the formal framework developed in \S \ref{s:reconprob}.  In particular, we explain conclusively how generalized sampling guarantees a stable and accurate reconstruction by deriving the exact values for $\mu$ and $\kappa$, and therefore $C$, as opposed to the nonsharp bounds given previously in \cite{BAACHShannon,BAACHAccRecov}.  Moreover, we reinterpret generalized sampling using geometry of Hilbert spaces, and in particular, the notions of oblique projections and subspace angles.  Next, we introduce a necessary and sufficient condition, the \textit{stable sampling rate}, which determines how to select the generalized sampling parameters so as to guarantee a good reconstruction.  This improves on the previous sufficient conditions of \cite{BAACHShannon,BAACHAccRecov}.  Once more, this condition is  easily computable, as explained in \S \ref{ss:stabsamp}.  We also discuss the connections between generalized sampling and computing with sections of infinite operators.

In \S \ref{s:optimality} we consider the question of optimality of generalized sampling.  That is, we pose the question: can another method outperform generalized sampling, and if so, in what sense?  Using the sharp bounds derived in \S \ref{s:GS}--\ref{ss:stabsamp}, we show that no method which is \textit{perfect} (a definition is given later) can exhibit better stability than generalized sampling.  Hence generalized sampling is an optimal, stable approach to the reconstruction problem amongst the class of perfect methods.  Moreover, for problems where the stable sampling rate grows linearly, we show that no method (perfect or nonperfect) can outperform generalized sampling in terms of the reconstruction accuracy by more than a constant factor.  Thus, although it is possible in theory to get a better approximation error with a different method, no method can converge at an asymptotically faster rate than generalized sampling.

In the final part of this paper, \S \ref{s:examples}, we consider the application of generalized sampling to the so-called \textit{uniform resampling problem}.  This problem concerns the computation of the Fourier coefficients of a function from nonuniformly-spaced samples of its Fourier transform.  We show that the standard approach to this problem is nothing more than an instance of consistent sampling, and we explain how in general this will lead to an exponentially large reconstruction constant $C$.  Next we consider the application of generalized sampling to this problem.  We prove that stable sampling rate is linear, and therefore generalized sampling is, in the senses defined in \S \ref{s:optimality}, an optimal, stable method for this problem.  Finally, we consider alternatives to uniform resampling, and show how the incorporation of different reconstruction systems -- specifically, splines and polynomials -- can lead to an improved reconstruction.

\subsection{Relation to sparsity and compressed sensing}
One of the most significant developments in signal and image reconstruction in the last several decades has been the introduction of sparsity-exploiting algorithms.  Techniques such as compressed sensing \cite{candesCSMag,donohoCS,EldarKutyniokCSBook,FornasierRauhutCS}, which exploit sparsity of the signal $f$ in a particular basis (wavelets, for example) to reduce the number of measurements required, have recently become extremely popular.  

Generalized sampling, in the form we discuss in this paper, does not exploit sparsity.  It guarantees recovery of all signals, sparse or otherwise, from sufficiently many of their measurements.  However, it transpires that generalized sampling can be combined with existing compressed sensing tools (randomization and convex optimization) to achieve subsampling, whenever the signal $f$ is sparse (or compressible) in the basis $\{ \phi_j \}_{j \in \bbN}$  \cite{BAACHGSCS}.   The importance of this development is that it allows for compressed sensing of analog signals (i.e.\ functions in function spaces) which have sparse (or compressible) information content in some infinite basis.  Conversely, the standard compressed sensing techniques and theorems apply only to finite-dimensional signals, i.e.\ vectors in finite-dimensional vector spaces.  We refer to \cite{BAACHGSCS} for details, and \cite{vetterliCSMag,VetterliFinInnov,VetterliFRI} for related methodologies based on analog, but finite information content, models for signals.

\section{The reconstruction problem}\label{s:reconprob}
We now describe the reconstruction problem in more detail.  To this end, suppose that $\{ \psi_{j} \}_{j \in \bbN}$ is a collection of elements of a separable Hilbert space $\rH$ (over $\bbC$) that forms a frame for a closed subspace $\rS$ of $\rH$ (the \textit{sampling} space).  In other words, $\spn \{ \psi_j : j \in \bbN \} $ is dense in $\rS$ and there exist constants $c_1,c_2>0$ (the frame constants) such that
\be{
\label{frameprop}
c_1 \| f \|^2 \leq \sum_{j \in \bbN} | \ip{f}{\psi_j} |^2 \leq c_2 \| f \|^2,\quad \forall f \in \rS,
}
where $\ip{\cdot}{\cdot}$ and $\nm{\cdot}$ are the inner product and norm on $\rH$ respectively \cite{christensen2003introduction}.  Suppose further that $\{ \phi_j \}_{j \in \bbN}$ is a collection of \textit{reconstruction} elements that form a frame for a closed subspace $\rT$ (the \textit{reconstruction} space), with frame constants $d_1,d_2>0$:
\be{
d_1 \| f \|^2 \leq \sum_{j \in \bbN} | \ip{f}{\phi_j} |^2 \leq d_2 \| f \|^2,\quad \forall f \in \rT.
}
Let $f \in \rH$ be a given element we wish to recover, and assume that we have access to the samples
\be{
\label{fsamples}
\hat{f}_{j} = \ip{f}{\psi_{j}},\quad j \in \bbN.
}
Note that the infinite vector $\hat{f} = \{ \hat{f}_j \}_{j \in \bbN}$ is an element of $\ell^2(\bbN)$.  Ignoring for the moment the issue of truncation -- namely, that in practice we only have access to the first $n$ measurements -- the reconstruction problem can now be stated as follows: given $\hat{f} = \{ \hat{f}_j \}_{j \in \bbN}$, find a reconstruction $\tilde{f}$ of $f$ from the subspace $\rT$.

\subsection{Stability and quasi-optimality}
There are two important conditions which a reconstruction, i.e.\ a mapping $\{ \hat{f}_j \}_{j \in \bbN} \mapsto \tilde{f}$, ought to possess.  The first is so-called \textit{quasi-optimality}:

\defn{
Let $F : \rH_0 \rightarrow \rT$, $f \mapsto \tilde f$ be a mapping, where $\rH_0$ is a subspace of $\rH$.  The quasi-optimality constant of $\mu = \mu(F) >0$ is the least number such that 
\bes{
\| f - \tilde f \| \leq \mu \| f - \cQ_{\rT} f \|,\quad \forall f \in \rH_0,
}
where $\cQ_{\rT}: \rH \rightarrow \rT$ is the orthogonal projection onto $\rT$.   If no such constant exists, we write $\mu = \infty$.  We say that $F$ is quasi-optimal if $\mu(F)$ is small.
}
Note that $\cQ_{\rT} f$ is the best, i.e.\ energy-minimizing, approximation to $f$ from $\rT$.  Thus, quasi-optimality states that the error committed by $\tilde f$ is within a small, constant factor of that of the energy-minimizing approximation.  The desire for quasi-optimal mappings arises from the fact that typical images and signals are known to be well represented in certain bases and frames, e.g.\ wavelets, splines or polynomials \cite{unser2000sampling}.  In other words, the error $\| f - \cQ_{\rT} f \|$ is small.  When reconstructing $f$ in the corresponding subspace $\rT$ from its measurements $\{ \hat{f}_j \}_{j \in \bbN}$ it is therefore vital that $\mu \ll \infty$.  Otherwise, the beneficial property of $\rT$ for the signal $f$ may be lost when computing the reconstruction $\tilde f$.

The second important consideration is that of stability.  For this, we introduce a condition number:
\defn{
\label{d:condnumb}
Let $\rH_0$ be a closed subspace of $\rH$ and suppose that $F : \rH_0 \rightarrow \rH$ is a mapping such that, for each $f \in \rH_0$, $F(f)$ depends only on the vector of samples $\hat{f} \in \ell^2(\bbN)$.  The (absolute) condition number $\kappa = \kappa(F)$ is given by
\be{
\label{condnumb}
\kappa = \sup_{f \in \rH_0} \lim_{\epsilon \rightarrow 0^+} \sup_{\substack{g \in \rH_0 \\ 0 < \| \hat{g} \|_{\ell^2} \leq \epsilon}} \left \{ \frac{\| F(f+g)-F(f) \|}{\| \hat{g} \|_{\ell^2}} \right \}.
}
We say that the mapping $F$ is well-conditioned if $\kappa$ is small.  Otherwise it is ill-conditioned.
}
A well-conditioned mapping $F$ is robust towards perturbations in the inputs $\{\hat{f}_j \}_{j \in \bbN}$.  Thus, in practice, where one always deals with noisy data, it is vitally important to have such a property.

It is worth noting at this stage that the condition number \R{condnumb} does not assume linearity of the mapping $F$.  If this is the case, then one has the much simpler form
\bes{
\kappa(F) = \sup_{\substack{f \in \rH_0 \\ \hat{f} \neq 0}} \left \{ \frac{\| F(f) \|}{\| \hat{f} \|} \right \}.
}
We also remark that \R{condnumb} is the \textit{absolute} condition number, as opposed to the somewhat more standard \textit{relative} condition number \cite{TrefethenBau}.  This is primarily for simplicity in the presentation: under some assumptions, it is possible to adapt the results we prove later in this paper for the latter.  

Occasionally, we will also consider the absolute condition number \textit{at} an element $f \in \rH_0$:
\be{
\label{local_condition}
\kappa_{f}(F) = \lim_{\epsilon \rightarrow 0^+} \sup_{\substack{g \in \rH_0 \\ 0 < \| \hat{g} \|_{\ell^2} \leq \epsilon}} \left \{ \frac{\| F(f+g)-F(f) \|}{\| \hat{g} \|} \right \},\quad f \in \rH_0.
}
This measures the local conditioning of $F$ around $f$.  Naturally, one has $\kappa(F) = \sup_{f \in \rH_0} \kappa_f(F)$.

For convenience, it is useful to introduce the notion of a reconstruction constant for $F$:
\defn{
Let $F$ be as in Definition \ref{d:condnumb}.  The reconstruction constant $C = C(F)$ is defined by
\bes{
C(F) = \max \left \{ \kappa(F) , \mu(F) \right \},
}
where the quantities $\mu(F)$ and $\kappa(F)$ are the quasi-optimality constant and condition number of $F$ respectively.  If $F$ is not quasi-optimal or if $\kappa(F)$ is not defined, then we set $C(F) = \infty$.
}

\subsection{The computational reconstruction problem}\label{ss:compreconprob}
As mentioned, in practice we do not have access to the infinite vector of samples $\hat{f}$.  Thus, the computational reconstruction problem concerns the recovery of $f$ from its first $n$ measurements $\hat{f}_1,\ldots,\hat{f}_n$.  Since we only have access to these samples, it is natural to consider finite-dimensional subspaces of $\rT$.  Thus, we let $\{ \rT_n \}_{n \in \bbN}$ be a sequence of finite-dimensional subspaces satisfying
\be{
\label{Tcond1}
\rT_n \subseteq \rT,\quad \mbox{$\dim(\rT_n) < \infty$},
}
and 
\be{
\label{Tcond2}
\cQ_{\rT_n} \rightarrow \cQ_{\rT},\quad n \rightarrow \infty,
}
strongly on $\rH$.  In other words, the spaces $\{ \rT_n \}_{n \in \bbN}$ form a sequence of finite-dimensional approximations to $\rT$.  Strictly speaking, the second condition is not necessary.  However, it is natural make this assumption in order to guarantee a good approximation.

\rem{
It is quite common in practice to define $\rT_n = \spn \{ \phi_1,\ldots,\phi_n \}$ to be the space spanned by the first $n$ elements of an infinite frame $\{ \phi_j \}_{j \in \bbN}$ for $\rT$.  Note that \R{Tcond1} and \R{Tcond2} automatically hold in this case.  Moreover, one has the nesting property $\rT_1 \subseteq \rT_{2} \subseteq \ldots$.  However, this is not necessary.  The reconstructions we consider in this paper are actually independent of the spanning system for $\rT_n$.  Such a system only needs to be specified in order to perform computations.  Hence, we consider the more general setting outlined above.
}

With this in hand, the computational reconstruction problem is now as follows: given the samples $\hat{f}_1,\ldots,\hat{f}_n$, compute a reconstruction to $f$ from the subspace $\rT_n$.  

When considering methods $F_n$ for this problem, it is clear that the constants reconstruction $C(F_n)$ should not grow rapidly with $n$.  If this is not the case, then increasing the number of measurements could, for example, lead to a worse approximation and increased sensitivity to noise.  We shall see examples of this in \S \ref{ss:consistfail}.  To avoid this scenario, we now make the following definition:

\defn{
\label{d:reconschemestab}
For each $n \in \bbN$, let $F_n$ be a mapping such that, for each $f$, $F_n(f)$ belongs to a finite-dimensional reconstruction space $\rT_n$ and depends only on the samples $\hat{f}^{[n]} = \{ \hat{f}_1,\ldots,\hat{f}_n \}$.  We say that the reconstruction scheme $\{F_n \}_{n \in \bbN}$ is numerically stable and quasi-optimal if  
\bes{
C^* : = \sup_{n \in \bbN} C(F_n)  < \infty,
}
where $C(F_n)$ is the reconstruction constant of $F_n$.  We refer to the constant $C^*$ as the reconstruction constant of the reconstruction scheme $\{ F_n \}_{n \in \bbN}$.
}

This definition incorporates the issue of \textit{approximation} into a sequence of reconstruction schemes.  Although in practice one only has access to a number of samples, it is natural to consider the behaviour of $F_n$ as $n$ -- the number of samples -- increases.  Ideally we would like $F_n(f)$ to behave like $\cQ_{n} f$, the best approximation to $f$ from $\rT_{n}$.  Namely, $F_n(f)$ should converge to $f$ at precise the same \textit{rate} as $\cQ_n f$.  This is vitally important from a practical standpoint.  The premise for computing a consistent reconstruction is the knowledge that $f$ is well represented in terms of the reconstruction system $\{\phi_j\}_{j \in \bbN}$.  This is equivalent to the property that the orthogonal projections $\cQ_n f$ converge rapidly.  Hence it is vital that the computed reconstruction $F_n(f)$ does not possess dramatically different behaviour as $n \rightarrow \infty$.  Put simply, there is little point reconstructing in the basis $\{ \phi_j \}_{j \in \bbN}$ if the good approximation properties of $f$ in this basis are destroyed by reconstruction technique.

\rem{
\label{r:relax}
In some applications, one may wish to relax the above definition slightly to allow mild growth of $C(F_n)$.  If $a_n \in (0,\infty)$ is an increasing sequence, we say that $\{ F_n \}_{n \in \bbN}$ is stable and quasi-optimal with respect to $\{ a_n \}_{n \in \bbN}$ if 
\bes{
C^* = \limsup_{n \rightarrow \infty} \frac{C(F_n)}{a_n} < \infty.
}
In other words, $C(F_n)$ can grow at worst like $\ord{a_n}$ as $n \rightarrow \infty$.
}

\section{Consistent reconstructions and oblique projections}\label{s:consistent}
We now consider the consistent sampling technique of \cite{eldar2003FAA,eldar2003sampling,eldar2005general,unser1994general,unserzerubia}.

\subsection{Consistent sampling}
Let us first return to the problem of recovering $f$ from its infinite vector of samples $\hat{f}$.  A simple and elegant way to obtain a reconstruction $F$ with small constant $C(F)$ is by solving the so-called \textit{consistency} conditions.  Specifically, we define $F(f) = \tilde f$ by 
\be{
\label{consistconds}
\langle \tilde{f} , \psi_j \rangle = \langle f , \psi_j \rangle,\quad j=1,2,\ldots ,\qquad \tilde f \in \rT.
}
Consistency means that the samples of $\tilde f$ agree with those of $f$.  We say that $\tilde f$ is a \textit{consistent reconstruction} of $f$, and refer to the mapping $F: f \mapsto \tilde f$ as \textit{consistent sampling}.

An analysis of consistent reconstructions, which we shall recap and extend in \S \ref{ss:consistentanalysis}, was given in \cite{eldar2003FAA,eldar2003sampling,eldar2005general}.  Crucial to this is the notion of oblique projections in Hilbert spaces, which we discuss next.  This tool will also be used later in analysing the generalized sampling technique.

\subsection{Oblique projections and subspace angles}
We commence with the definition of a subspace angle:

\defn{
Let $\rU$ and $\rV$ be closed subspaces of a Hilbert space $\rH$ and $\cQ_{V} : \rH \rightarrow \rV$ the orthogonal projection onto $\rV$.  The subspace angle $\theta = \theta_{\rU \rV} \in [0,\frac{\pi}{2}]$ between $\rU$ and $\rV$ is given by
\be{
\label{subangle}
\cos (\theta_{\rU \rV}) = \inf_{\substack{u \in \rU \\ \| u \|=1}} \| \cQ_{V} u \|.
}
}
Note that there are a number of different ways to define the angle between subspaces \cite{steinberg2000oblique,Tang1999Oblique}.  However, \R{subangle} is the most convenient for this paper.  We shall also make use of the following equivalent expression for $\cos \left ( \theta_{\rU \rV} \right )$:
\be{
\label{subangle2}
\cos \left ( \theta_{\rU \rV} \right ) = \inf_{\substack{u \in \rU \\ \| u \| = 1}} \sup_{\substack{v \in \rV \\ \| v \| = 1}} \ip{u}{v}.
}
We are interested in subspaces for which the cosine of the associated angle is nonzero.  The following lemma is useful in this extent:
\lem{
\label{l:subspace_cond_equivalent}
Let $\rU$ and $\rV$ be closed subspaces of a Hilbert space $\rH$.  Then $\cos \left ( \theta_{\rU \rV^{\perp}} \right ) > 0$ if and only if $\rU \cap \rV = \{ 0 \}$ and $\rU + \rV$ is closed $\rH$.
}
\prf{
See \cite[Thm. 2.1]{Tang1999Oblique}.
}
We now make the following definition:
\defn{
Let $\rU$ and $\rV$ be closed subspaces of a Hilbert space $\rH$.  Then $\rU$ and $\rV$ satisfy the subspace condition if $\cos \left ( \theta_{\rU \rV^{\perp} } \right ) > 0$, or equivalently, if $\rU \cap \rV = \{ 0 \}$ and $\rU + \rV$ is closed in $\rH$.
}
Subspaces $\rU$ and $\rV$ satisfying this condition give rise to a decomposition $\rU \oplus \rV = \rH_0$ of a closed subspace $\rH_0$ of $\rH$.  Equivalently this ensures the existence of a projection of $\rH_0$ with range $\rU$ and kernel $\rV$.  We refer to such a projection as an \textit{oblique projection} and denote it by $\cW_{\rU \rV}$.  Note that $\cW_{\rU \rV}$ will not, in general, be defined over the whole of $\rH$, but rather the subspace $\rH_0$.  However, this is true whenever $\rV = \rU^{\perp}$, for example, and in this case the projection $\cW_{\rU \rV}$ coincides with the orthogonal projection $\cQ_{\rU}$.

We shall also require the following results on oblique projections (see \cite{BuckholtzIdempotents,SzyldOblProj}):

\thm{
\label{t:oblproj}
Let $\rU$ and $\rV$ be closed subspaces of $\rH$ with $\rU \oplus \rV = \rH$.  Then
\bes{
\| \cW_{\rU \rV} \| =  \| \cI - \cW_{\rU \rV} \| = \sec \left (  \theta_{\rU \rV^{\perp}} \right ), 
}
where $\nm{\cdot}$ is the standard norm on the space of bounded operators on $\rH$.
}

\cor{
\label{c:consist}
Suppose that $\rU$ and $\rV$ are closed subspaces of $\rH$ satisfying the subspace condition, and let $\cW_{\rU \rV} : \rH_0 \rightarrow \rU$ be the oblique projection with range $\rU$ and kernel $\rV$, where $\rH_0 = \rU \oplus \rV$.  Then
\be{
\label{consiststab}
\| \cW_{\rU \rV} f \| \leq \sec \left ( \theta_{\rU \rV^\perp} \right ) \| f \|,\quad \forall f \in \rH_0,
}
and 
\be{
\label{consisterr}
\| f - \cQ_{\rU} f \| \leq \| f - \cW_{\rU \rV} f \| \leq \sec \left ( \theta_{\rU \rV^\perp} \right ) \| f -\cQ_{\rU} f\|,\quad \forall f \in \rH_0,
}
where $\cQ_{\rU} : \rH \rightarrow \rU$ is the orthogonal projection.  Moreover, the upper bounds in \R{consiststab} and \R{consisterr} are sharp.
}
\prf{
The sharp bound \R{consiststab} follows immediately from Theorem \ref{t:oblproj}.  For \R{consisterr} we first observe that $(\cI-\cW_{\rU \rV}) = (\cI-\cW_{\rU \rV})(\cI - \cQ_{\rU})$, since $\cW_{\rU \rV}$ and $\cQ_{\rU}$ are both projections onto $\rU$.  Hence, by Theorem \ref{t:oblproj},
\bes{
\| f - \cW_{\rU \rV} f \| = \| (\cI - \cW_{\rU \rV})(\cI - \cQ_{\rU}) f \| \leq \sec \left (\theta_{\rU \rV^{\perp}} \right ) \| f - \cQ_{\rU} f \|,
}
with sharp bound.
}

\rem{
Although arbitrary subspaces $\rU$ and $\rV$ need not obey the subspace condition, this is often the case in the important examples arising in practice.  For example, if $\rU \subseteq \rV^{\perp}$ then it follows immediately from the definition that $\cos \left ( \theta_{\rU \rV^{\perp}} \right ) =1$. 
}

To complete this section, we present the following lemma which will be useful in what follows:

\lem{
\label{l:findimoblique}
Let $\rU$ and $\rV$ be closed subspaces of $\rH$ satisfying the subspace condition.  Suppose also that $\dim (\rU) = \dim (\rV^{\perp}) =n < \infty$.  Then $\rU \oplus \rV = \rH$.
}
\prf{
Note that $\rU \oplus \rV = \rH$ if and only if $\cos \left ( \theta_{\rU \rV^{\perp} } \right)$ and $\cos \left ( \theta_{\rV^{\perp} \rU} \right )$ are both positive \cite[Thm.\ 2.3]{Tang1999Oblique}.  Since $\cos \left ( \theta_{\rU \rV^{\perp} } \right) > 0$ by assumption, it remains to show that $\cos  \left ( \theta_{\rV^{\perp} \rU} \right ) > 0$.  Consider the mapping $\cQ_{\rV^{\perp}} \big |_{\rU} : \rU \rightarrow \rV^{\perp}$.  We claim that this mapping is invertible.  Since $\rU$ and $\rV^{\perp}$ have the same dimension it suffices to show that $\cQ_{\rV^{\perp}} \big |_{\rU}$ has trivial kernel.  However, the existence of a nonzero $u \in \rU$ with $\cQ_{\rV^{\perp}} u = 0$ implies that $\cos  \left ( \theta_{\rU \rV^{\perp}} \right )=0$; a contradiction.  Thus $\cQ_{\rV^{\perp}} \big |_{\rU}$ is invertible, and in particular, it has range $\rV^{\perp}$.  Now consider $\cos  \left ( \theta_{\rV^{\perp} \rU} \right )$.  By \R{subangle2} and this result,
\eas{
\cos  \left ( \theta_{\rV^{\perp} \rU} \right ) &= \inf_{\substack{w \in \rV^{\perp} \\ w \neq 0}} \sup_{\substack{u \in \rU \\ u \neq 0}} \frac{\ip{w}{u}}{\| w \| \| u \|} 
= \inf_{\substack{u' \in \rU \\ u' \neq 0}} \sup_{\substack{u \in \rU \\ u \neq 0}} \frac{\ip{\cQ_{\rV^{\perp}} u'}{u}}{\| \cQ_{\rV^{\perp}} u' \| \| u \|} 
 \geq \inf_{\substack{u' \in \rU \\ u' \neq 0}} \frac{\| \cQ_{\rV^{\perp}} u' \|}{\| u' \|}
 = \cos \left ( \theta_{\rU \rV^{\perp}} \right ) > 0.
}  
This completes the proof.
}

\subsection{Quasi-optimality of consistent sampling}\label{ss:consistentanalysis}

Oblique projections arise in many types of sampling \cite{aldroubi1996oblique}, and consistent reconstructions are intimately related with such mappings.  The following result was proved in \cite[Thm.\ 2.1]{eldar2005general}:

\lem{
\label{l:consistoblique}
Suppose that $\rT$ and $\rS^{\perp}$ satisfy the subspace condition.  If $f \in \rH_0 = \rT \oplus \rS^{\perp}$ then there exists a unique $\tilde{f} \in \rT$ satisfying \R{consistconds}.  Specifically, the mapping $F: \rH_0 \rightarrow \rT, f \mapsto \tilde{f}$ coincides with the oblique projection $\cW_{\rT\rS^{\perp}} $.
}
As a result of this lemma, consistent reconstructions are equivalent to oblique projections.  However, we note one important distinction.  When defining the consistent reconstruction \R{consistconds}, we assume that a frame $\{ \psi_j \}^{\infty}_{j=1}$ of $\rS$ is given.  Indeed, this is natural in view of the sampling process.  However, the oblique projection $\cW_{\rT\rS^{\perp}}$, being determined solely by the spaces $\rT$ and $\rS$, is actually independent of this basis.  In light of Lemma \ref{l:consistoblique}, the same must also be true for $\tilde f$.  In fact, as we detail in \S \ref{ss:consistcompute}, specification of frames or bases $\{ \phi_j \}^{\infty}_{j=1}$ and $\{ \psi_j \}^{\infty}_{j=1}$ is only necessary when writing \R{consistconds} as a linear system of equations to be solved numerically.

Lemma \ref{l:consistoblique}, in combination with Corollary \ref{c:consist}, gives the following sharp bounds for consistent sampling:
\ea{
\| \tilde{f} \| &\leq \sec \left ( \theta_{\rT \rS } \right ) \| f \|,\quad \forall f \in \rH_0,\label{consiststabA}
\\
\| f - \cQ_{\rT} f \| \leq \| f - \tilde{f} \| &\leq \sec \left ( \theta_{\rT \rS} \right ) \| f -\cQ_{\rT} f\|,\quad \forall f \in \rH_0.\label{consisterrA}
}
The latter illustrates quasi-optimality of the reconstruction, whereas the former gives a continuous stability estimate for $\tilde f$, i.e.\ the norm of the reconstruction is bounded by a constant multiple of the norm of the input signal.  Note that \R{consiststabA} and \R{consisterrA} were derived previously in \cite{unserzerubia} and \cite{unser1994general} respectively.  However, the observation that they are also sharp does not, to the best of our knowledge, appear in the literature on consistent reconstructions.

Another property of the consistent reconstruction is confirmed by the above bounds.  Namely, it is a \textit{perfect} reconstruction:

\defn{
\label{d:perfect}
A mapping $F : \rH \rightarrow \rT$ is a perfect reconstruction if (i) for each $f \in \rH$, $F(f)$ depends only on the vector of samples $\hat{f}$,  and (ii) $F(f) = f$ whenever $f \in \rT$.
}

\subsection{The condition number of consistent sampling}
The bound \R{consisterrA} demonstrates that the quasi-optimality constant of consistent sampling is $\mu(F) = \sec \left ( \theta_{\rT \rS} \right )$.  We now wish to determine the condition number.  For this, it is useful to first recall several basic facts about frames \cite{christensen2003introduction}.  Given the sampling frame $\{ \psi_j \}_{j \in \bbN}$ for the subspace $\rS$, we define the \textit{synthesis} operator $\cS : \ell^2(\bbN) \rightarrow \rH$ by
\bes{
\cS \alpha= \sum_{j \in \bbN} \alpha_j \psi_j,\quad \alpha = \{ \alpha_j \}_{j \in \bbN}  \in \ell^2(\bbN).
}
Its adjoint, the \textit{analysis} operator, is defined by
\bes{
\cS^* f = \hat{f} = \{ \ip{f}{\psi_j} \}_{j \in \bbN},\quad f \in \rH.
}
The resulting composition $\cP = \cS \cS^* : \rH \rightarrow \rH$, given by
\be{
\label{Pdef}
\cP f = \sum_{j \in \bbN} \ip{f}{\psi_j} \psi_j,\quad \forall f \in \rH,
}
is well-defined, linear, self-adjoint and bounded.  Moreover, the restriction $\cP |_{\rS} : \rS \rightarrow \rS$ is positive and invertible with $c_1 \cI |_{\rS} \leq \cP |_{\rS} \leq c_2 \cI |_{\rS}$.  We now require the following lemma:
\lem{
\label{l:Pbdd}  Suppose that $\rT$ and $\rS^{\perp}$ satisfy the subspace condition, and let $\cP$ be given by \R{Pdef}.  Then
\be{
\label{E:Pbdd}
 c_1 \cos^2 (\theta_{\rT \rS}) \hspace{1mm} \cI |_{\rT} \leq  \cP |_{\rT} \leq c_2  \cI |_{\rT},
}
where $c_1,c_2$ are the frame constants appearing in \R{frameprop}.
}
\prf{
Let $f \in \rH$ be arbitrary, and write $f = \cQ_{\rS} f + \cQ_{\rS^\perp} f$.  Then
\be{
\label{Pdecomp}
\ip{\cP f}{f} = \sum_{j \in \bbN} | \ip{f}{\psi_j} |^2 = \sum_{j \in \bbN} | \ip{\cQ_{\rS} f}{\psi_j} |^2 = \ip{\cP \cQ_{\rS} f}{\cQ_{\rS} f }.
}
Suppose now that $\phi \in \rT$.  Using \R{Pdecomp} and the frame condition \R{frameprop} we find that
\bes{
c_1 \| \cQ_{\rS} \phi \|^2  \leq \ip{\cP \phi}{\phi} \leq c_2 \| \cQ_{\rS} \phi \|^2.
}
The upper bound in \R{E:Pbdd} follows immediately from the observation  that $\| \cQ_{\rS} \phi \| \leq \| \phi \|$.  For the lower bound we use the definition of the subspace angle $\theta_{\rT \rS}$, and the fact that $\phi \in \rT$.
}

\cor{
\label{c:consistcondition} Suppose that $\rS$ and $\rT$ are as in Lemma \ref{l:Pbdd}, and let $F : \rH_0 := \rT \oplus \rS^{\perp} \rightarrow \rT$ denote the consistent reconstruction defined by \R{consistconds}.  Then the condition number $\kappa(F)$ satisfies
\bes{
\frac{\sec (\theta_{\rT \rS})}{\sqrt{c_2}}   \leq \kappa(F) \leq \frac{ \sec  (\theta_{\rT \rS})}{\sqrt{c_1}}.
}
}
\prf{
Since the reconstruction $\tilde{f} = F(f) \in \rT$ is defined by \R{consistconds}, we have
\bes{
\| \hat{f} \|^2_{\ell^2} = \sum_{j \in \bbN} | \ip{\tilde{f}}{\psi_j} |^2 =\ip{\cP \tilde{f}}{\tilde{f}}.
}
Hence, by the previous lemma, $\| \hat{f} \|^2_{\ell^2}  \geq c_1 \cos^2 (\theta_{\rT \rS}) \| \tilde{f} \|^2$.  Since $F$ is linear, this now gives
\bes{
\kappa(F) = \sup_{\substack{f \in \rH_0 \\ \hat{f} \neq0}} \left \{ \frac{\| F(f) \|}{\|\hat{f} \|_{\ell^2}} \right \} \leq \frac{\sec (\theta_{\rT \rS})}{\sqrt{c_1}}. 
}
On the other hand, since the reconstruction $F$ is perfect,
\bes{
\kappa(F) \geq \sup_{\substack{f \in \rT \\ \hat{f} \neq0}}\left \{ \frac{\| f \|}{\|\hat{f} \|_{\ell^2}} \right \} = \sup_{\substack{f \in \rT \\ f \neq0}}\left \{ \frac{\| f \|}{\|\hat{f} \|_{\ell^2}} \right \}.
}
Moreover, by \R{Pdecomp}, we have $\| \hat{f} \|^2_{\ell^2} \leq c_2 \| \cQ_{\rS} f \|^2$.  Hence
\bes{
\kappa(F) \geq \frac{1}{\sqrt{c_2}} \sup_{\substack{f \in \rT \\ f \neq0}} \left \{  \frac{\| f \|}{\| \cQ_{\rS} f \|}  \right \} = \frac{\sec (\theta_{\rT \rS})}{\sqrt{c_2}},
}
as required.
}

Combining Corollaries \ref{c:consist} and \ref{c:consistcondition}, we now find that the reconstruction constant $C(F)$ of consistent sampling satisfies
\be{
\label{inf_consist_recon_const}
\sec (\theta_{\rT \rS}) \max \{ 1,\tfrac{1}{\sqrt{c_2}} \} \leq C(F) \leq \sec (\theta_{\rT \rS}) \max \{ 1 , \tfrac{1}{\sqrt{c_1}} \}.
}
Hence, if $\rS$ and $\rT$ are not close to perpendicular (i.e.\ if $\cos (\theta_{\rT \rS})$ is not too small), then consistent sampling is stable and quasi-optimal.

\subsection{Consistent sampling for the computational reconstruction problem}\label{ss:consistcompute}
The consistent reconstruction $\tilde{f}$ solves the reconstruction problem of recovering $f$ from the infinite vector $\hat{f} = \{ \hat{f}_j \}_{j \in \bbN}$.  Note that once a frame $\{ \phi_j \}_{j \in \bbN}$ is specified for $\rT$, this is equivalent to the infinite system of linear equations $U \alpha = \hat{f}$, where $U$ is the infinite matrix
\be{
\label{Udef}
U = 
  \left(\begin{array}{ccc} \left < \phi_1 , \psi_1 \right >   & \left < \phi_2 , \psi_1 \right > & \cdots \\  
\left < \phi_1 , \psi_2 \right >   & \left < \phi_2 , \psi_2\right > & \cdots \\ 
\vdots  & \vdots  & \ddots   \end{array}\right),
}
and $\alpha = \{ \alpha_j \}_{j \in \bbN} \in \ell^2(\bbN)$ is such that $\tilde{f} = \sum_{j \in \bbN} \alpha_j \phi_j$.  Observe that $U$, which we may view as an operator on $\ell^2(\bbN)$, coincides with $\cX = \cS^* \cT : \ell^2(\bbN) \rightarrow \ell^2(\bbN)$, where $\cS$ and $\cT$ are the synthesis operators for the frames $\{ \psi_j \}^{\infty}_{j=1}$ and $\{ \phi_j \}^{\infty}_{j=1}$ respectively.

Clearly this approach does not solve the computational reconstruction problem outlined in \S \ref{ss:compreconprob} since one cannot compute solutions to $U \alpha = \hat{f}$ in general.  To overcome this, the  standard approach \cite{eldar2003FAA,eldar2003sampling,EldarMinimax,UnserHirabayashiConsist,unser2000sampling} is to replace the infinite consistency conditions \R{consistconds} by a finite version.  That is, we seek a reconstruction $\tilde{f}_{n,n}$ defined by
 \be{
\label{E:findimconsist}
\ip{\tilde f_{n,n}}{\psi_j} = \ip{f}{\psi_j},\quad j=1,\ldots,n,\qquad \tilde f_{n,n} \in \rT_n,
}
(the use of the double index in $\tilde{f}_{n,n}$ is for agreement with subsequent notation).  Note that if $\rT_n = \spn \{ \phi_1,\ldots,\phi_n \}$ then this is equivalent to the finite-dimensional linear system $U^{[n,n]} \alpha^{[n,n]} = \hat{f}^{[n]}$, where
\be{
\label{Umdef}
U^{[n,n]} = 
  \left(\begin{array}{ccc} \left < \phi_1 , \psi_1 \right >    & \cdots & \left < \phi_n , \psi_1 \right > \\  
\vdots  &  \ddots & \vdots  \\
\left < \phi_1 , \psi_n \right >    & \cdots & \left < \phi_n , \psi_n\right >  \end{array}\right),
}
$\hat{f}^{[n]} = \{ \hat{f}_1,\ldots,\hat{f}_n \} $, $\alpha^{[n,n]} = \{ \alpha^{[n]}_1 ,\ldots,\alpha^{[n]}_n \}$ and $\tilde{f}_{n,n}$ is given by $\sum^{n}_{j=1} \alpha^{[n,n]}_{j} \phi_j$.

The condition \R{E:findimconsist} is completely natural and reasonable to enforce: it states that the reconstruction $\tilde{f}_{n,n}$ is consistent with the available data $\{ \hat{f}_j \}^{n}_{j=1}$.  Moreover, it is tempting to think that stability and quasi-optimality of the infinite-dimensional consistent reconstruction $\tilde{f}$ should imply the same behaviour of $\tilde{f}_{n,n}$.  In other words, $C^* = \sup_{n \in \bbN} C(F_{n,n})$ should be both finite and not too large.  However, as we explain in the next section, there is no guarantee that this will be the case in practice.

First, however, we determine the reconstruction constant for this approach.  As a direct consequence of Lemmas \ref{l:findimoblique}, \ref{l:consistoblique} and Corollary \ref{c:consist}, we have

\cor{
\label{c:finconsist} Let $\rS_n = \spn \{ \psi_1,\ldots,\psi_n \}$ and suppose that
\be{
\label{finspacesass}
\cos \left ( \theta_{n,n} \right ) > 0,
}
where $\theta_{n,n} = \theta_{\rT_n \rS_n}$.  Then, for each $f \in \rH_n : = \rT_n \oplus \rS^{\perp}_n$ there exists a unique $\tilde f_{n,n} \in \rT_n$ satisfying \R{E:findimconsist}.  Moreover, the mapping $F_{n,n}: \rH_n \rightarrow \rT_n, f  \mapsto \tilde{f}_{n,n}$ coincides with the oblique projection $\cW_{\rT_n  \rS^{\perp}_n}$, and we have the sharp bounds
\be{
\label{finstab}
\| \tilde f_{n,n} \| \leq \sec \left ( \theta_{n,n} \right ) \| f \|,
}
and
\be{
\label{finerrest}
\| f - \cQ_{n} f \| \leq \| f - \tilde f_{n} \| \leq \sec \left ( \theta_{n,n} \right ) \| f - \cQ_{n} f\|,
}
where $\cQ_{n}$ is the orthogonal projection onto $\rT_{n}$.  If $\dim(\rT_n) = \dim(\rS_n)$ (in particular, if both $\{ \psi_j \}_{j \in \bbN}$ and $\{ \phi_j \}_{j \in \bbN}$ are bases), then the above conclusions hold with $\rH_0 =  \rH$.
}
This theorem demonstrates that the quasi-optimality constant of consistent sampling satisfies
\be{
\label{consistmu}
\mu(F_{n,n}) = \sec \left ( \theta_{n,n} \right ).
}
We now state the following result concerning the condition number (a proof is given in \S \ref{s:GS}):
\cor{
\label{c:reconconst}
Suppose that $\cos (\theta_{n,n}) > 0$ and that $\dim (\rS_n) = \dim(\rT_n)$.  Then the condition number of consistent sampling satisfies
\bes{
\kappa(F_{n,n}) \geq \frac{1}{\sqrt{c_2}} \sec (\theta_{n,n}).
}
}
This result, in combination with \R{consistmu}, implies that the reconstruction constant $C(F_{n,n})$ satisfies
\bes{
C(F_{n,n}) \geq \max \{ 1 , \tfrac{1}{\sqrt{c_2}} \} \sec ( \theta_{n,n} ).
}
Hence, good behaviour of the reconstructions (i.e.\ stability and quasi-optimality for all $n$) is only guaranteed if the subspace angles $\theta_{n,n}$ remain bounded away from $\frac{\pi}{2}$ for all $n \in \bbN$.  As we next explain, there is no need for this to be the case in general, even when the angle $\theta_{\rT \rS} $ between the infinite-dimensional subspaces $\rT$ and $\rS$ satisfies this condition.

\subsection{Example}\label{ss:consistfail}

It is easiest to illustrate this issue by introducing the main example we shall consider in this paper: namely, the \textit{uniform resampling} (URS) problem.  This problem, which we shall discuss in further detail in \S \ref{s:examples}, occurs in the reconstruction of signals and images from nonuniformly sampled Fourier data.

An example of this problem, written in terms of the consistent sampling framework, is as follows.  
Let $\rH = \rL^2(-1,1)$ and, for $0 < \delta < 1$ set
\bes{
\psi_{j}(x) = \frac{1}{\sqrt{2}}\E^{\I j \pi \delta x},\quad \phi_j(x) = \frac{1}{\sqrt{2}} \E^{\I j \pi x},\quad j \in \bbZ,
}
(for convenience we now index over $\bbZ$, as opposed to $\bbN$).  In this case, $\{ \psi_j \}_{j \in \bbZ}$ is a tight frame for $\rS = \rH$ with frame bounds $c_1 =c_2 = \frac{1}{\delta}$, and $\{ \phi_j \}_{j \in \bbZ}$ is an orthonormal basis for $\rT = \rH$.  Since $\rS = \rT = \rH$ we have $\theta_{\rT \rS} = 0$.

\begin{figure}
\begin{center}
$\begin{array}{ccc}
 \includegraphics[width=4.75cm]{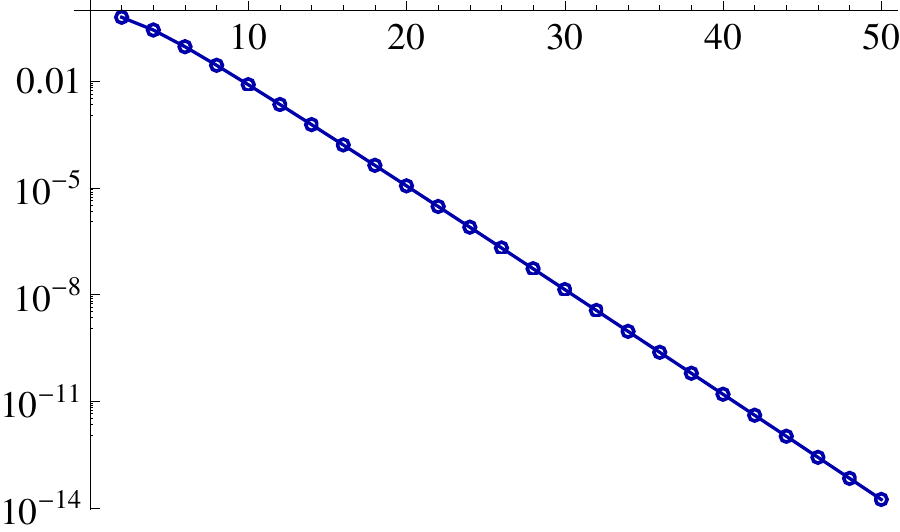}   &   \includegraphics[width=4.75cm]{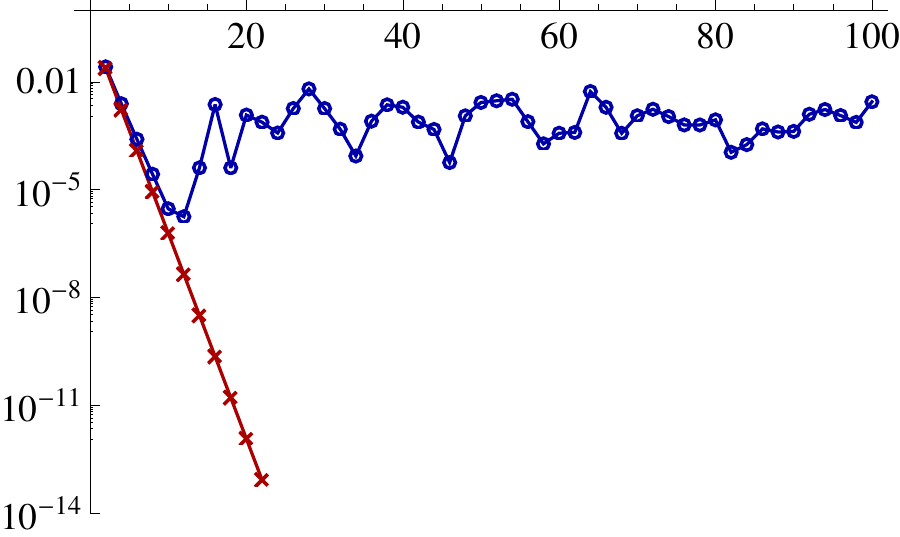} &   \includegraphics[width=4.75cm]{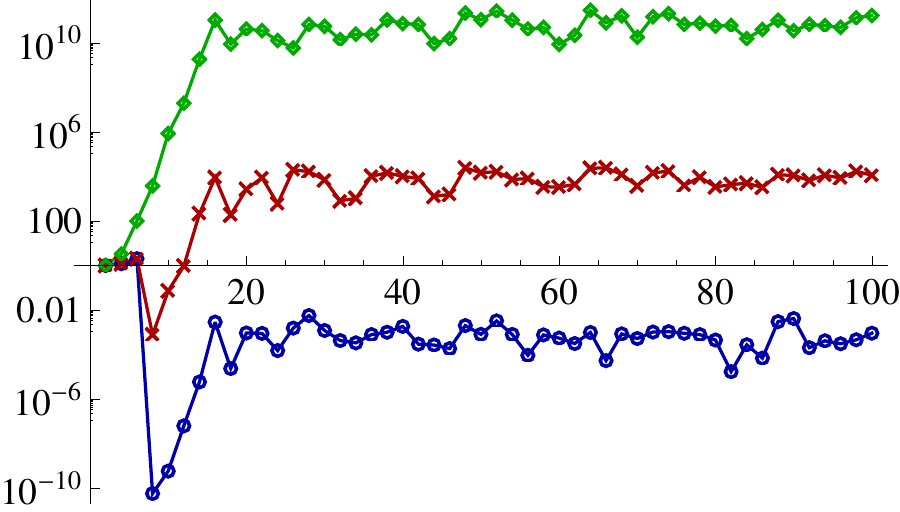} \\
 (a) & (b) & (c)
\end{array}$
\caption{\small (a): $\cos(\theta_{n,n})$ against $n$.  (b): the errors $\| f - \tilde{f}_{n,n} \|$ (circles) and $\| f - \cQ_{n} f \|$ (crosses) against $n$, where $\tilde{f}_{n,n}$ is the consistent reconstruction of $f$, $\cQ_{n} f$ is the orthogonal projection of $f$ onto $\rT_n$, and $f(x) = \frac{1}{2+\cos \pi x}$.  (c): the error $\| f - \tilde{f}_{n,n} \|$ against $n$, where $f(x) = \frac{1}{\sqrt{2}} \E^{8 \I \pi x}$ and each sample $\hat{f}_j$ is perturbed by an amount $\eta_j$ chosen uniformly at random with $| \eta_j | \leq \eta$, and $\eta = 0,10^{-9},10^{-2}$ (circles, crosses and diamonds respectively).}  \label{f:Ex1}
\end{center}
\end{figure}

In Figure \ref{f:Ex1}(a) we plot the behaviour of $\cos(\theta_{n,n})$ against $n$ for $\delta = \frac12$, where
\bes{
\rS_n = \spn \{ \psi_j : |j| \leq n \},\quad \rT_n = \spn \{ \phi_j : | j |\leq n \}.
}
As is evident, the quantity $\cos \left( \theta_{n,n} \right )$ is exponentially small in $n$.  In particular, when $n=50$ Figure \ref{f:Ex1}(a), in combination with Corollary \ref{c:reconconst}, implies that the reconstruction constant for the consistent reconstruction based on these spaces is around $10^{14}$ in magnitude.  Hence we expect extreme instability.

As mentioned, this example is an instance of the URS problem.  The reconstruction space $\rT_n$ is the space of trigonometric polynomials of degree $n$, and as such, it is particularly well-suited for recovering smooth and periodic functions.  In particular, if $f$ is smooth, then the error $\| f - \cQ_{n} f \|$ decays extremely rapidly: namely, faster than any algebraic power of $n^{-1}$.  However, such good approximation properties are destroyed when moving to $\tilde{f}_{n,n}$.  In Figure \ref{f:Ex1}(b) we plot the error $\| f - \tilde{f}_{n,n} \|$ for the consistent reconstruction, as well as that of the orthogonal projection $\cQ_n f$.  The latter decays rapidly with $n$, as expected.  On the other hand, the maximal achievable accuracy of the consistent reconstruction is limited to only around one or two digits due to the ill-conditioning and the effect of round-off errors.

The situation worsens significantly when the samples $\hat{f}_j$ are corrupted by noise.  The function used in Figure \ref{f:Ex1}(c) actually lies in $\rT_n$ whenever $n \geq 8$, and therefore, in theory at least, should be recovered perfectly.  However, this is completely obliterated by even moderate amounts of noise.  For example, even with noise at amplitude $10^{-9}$ the reconstruction error is around $10^4$ (i.e.\ an amplification of $10^{15}$), rendering such an approach useless for this particular problem.

\rem{
The exponential blow-up of the reconstruction constant in the above example is by no means unique to this particular problem.  In \cite{BAACHShannon} several other problems, based on spaces $\rT_n$ consisting of polynomials or piecewise constant functions, were shown to exhibit similar exponential growth.
}

\subsection{Operator-theoretic interpretation}\label{ss:optheory}
To sum up, the failure of the finite-dimensional consistent reconstruction $\tilde f_{n,n}$ is due to the poor behaviour of the finite subspace angles $\theta_{n,n} = \theta_{\rT_n,\rS_n}$ in relation to $\theta = \theta_{\rT \rS}$

Another interpretation of this failure  was provided in \cite{BAACHShannon}.  If the sampling and reconstruction bases are orthonormal (this fact is not necessary, but simplifies what follows) then $\cos (\theta_{n,n})$ and $\cos (\theta)$ coincide with the minimal singular values of the matrices $U^{[n,n]}$ and $U$ respectively.  Hence, the fact that $\theta_{n,n}$ may behave wildly, even when $\theta$ is bounded away from $\frac{\pi}{2}$, demonstrates that the spectra of the matrices $U^{[n,n]}$ poorly approximate the spectrum of $U$.  

This question -- namely, how well does a sequence of finite-rank operators approximate the spectrum of a given infinite-rank operator -- is one of the most fundamental in the field of spectral theory.  Recalling the definition \R{Umdef}, we notice that $U^{[n,n]}$ is nothing more than the $n \times n$ \textit{finite section} of $U$.  In other words, if $\{ e_j \}$ is the canonical basis for $\ell^2(\bbN)$ and $P_n : \ell^2(\bbN) \rightarrow \spn \{ e_1,\ldots,e_n \}$ is the orthogonal projection, then $U^{[n,n]} = P_n U P_n$.  Moreover, $\hat{f}^{[n]} = P_{n} \hat{f}$, and thus the finite-dimensional consistent reconstruction $\tilde f_{n,n}$ is precisely the result of the \textit{finite section method} applied to the equations $U \alpha = \hat{f}$.

The failure of consistent reconstructions can consequently be viewed from this perspective.  The properties of finite sections have been extensively studied over the last several decades \cite{bottcher1996,hansen2008,lindner2006}, and unfortunately there is no general guarantee they are well behaved.  To put this in a formal perspective, suppose for the moment that we approximate the operator $U$ with a sequence $U^{[n]}$ of finite-rank operators (which may or may not be finite sections), and instead of solving $U \alpha = \hat{f}$, we solve $U^{[n]} \alpha^{[n]} = \hat{f}^{[n]}$.  For obvious reasons, it is vitally important that this sequence satisfies the three following conditions:

\enum{
\item[(i)] \textit{Invertibility:} $U^{[n]}$ is invertible for all $n=1,2,\ldots$.
\item[(ii)] \textit{Stability:} $\| (U^{[n]})^{-1} \|$ is uniformly bounded for all $n=1,2,\ldots$.
\item[(iii)] \textit{Convergence:} the solutions $\alpha^{[n]} \rightarrow \alpha$ as $n \rightarrow \infty$.
}
Unfortunately, there is no guarantee that finite sections, and therefore the consistent reconstruction technique, possess any of these properties.  In fact, one requires rather restrictive conditions, such as positive self-adjointness, for this to be the case.  Typically operators $U$ of the form \R{Udef} are not self-adjoint, thereby making this approach unsuitable in general for discretizing $U \alpha = \hat{f}$. 

The framework of generalized sampling, which we next discuss, overcomes these issues by obtaining a sequence of operators that possess the properties (i)--(iii) above.  The key to doing this is to allow the number of samples $m$ to vary independently from the number of reconstruction vectors $n$.  When done in this way, we obtain a finite-dimensional operator $U^{[n,m]}$ (which now depends on both $m$ and $n$) that inherits the spectral structure of its infinite-dimensional counterpart $U$, provided $m$ is sufficiently large for a given $n$.  It turn, this ensures a stable, quasi-optimal reconstruction.

\section{Generalized sampling}\label{s:GS}
We now consider generalized sampling.  This framework was first introduced in \cite{BAACHShannon}, and applied to the resolution of the Gibbs phenomenon in \cite{BAACHAccRecov}.  Several extensions have also been pursued, to infinite-dimensional compressed sensing \cite{BAACHGSCS}, inverse and ill-posed problem \cite{AHHTillposed}, and problems where the sampling and reconstruction systems lie in different Hilbert space \cite{AdcockHansenSpecData}.

From now on we shall assume that the subspaces $\rT$ and $\rS^{\perp}$ obey the subspace condition.  In other words, $\cos \theta_{\rT \rS} > 0$.  Without this, the infinite-dimensional reconstruction problem is itself ill-posed, and thus it becomes far more difficult to obtain stable, quasi-optimal reconstructions.

Let $\rS_{m} = \spn \{ \psi_1,\ldots,\psi_m \}$ and suppose that $\{ \rT_n \}_{n \in \bbN}$ is a sequence of finite-dimensional reconstruction spaces satisfying \R{Tcond1} and \R{Tcond2}.  We seek a reconstruction $\tilde{f}_{n,m} \in \rT_n$ of $f$ from the $m$ samples $\hat{f}_1,\ldots,\hat{f}_m$.  Let $\cP_{m} : \rH \rightarrow \rS_{m}$ be the finite rank operator given by
\bes{
\cP_{m} g = \sum^{m}_{j=1} \ip{g}{\psi_j} \psi_j.
}
Note that, due to \R{frameprop}, the sequence of operators $\cP_{m}$ converge strongly to $\cP$ on $\rH$ \cite{christensen2003introduction}, where $\cP$ is given by (\ref{Pdef}).  With this to hand, the approach proposed in \cite{BAACHShannon} is to define $\tilde{f}_{n,m} \in \rT_{n}$ as the solution of the equations
\be{
\label{redconsisteqns}
\ip{\cP_{m} \tilde{f}_{n,m}}{\phi_j} = \ip{\cP_{m} f }{\phi_j},\quad j=1,\ldots,n,\qquad \tilde f_{n,m} \in \rT_n.
}
We refer to the mapping $F_{n,m} : f \mapsto \tilde f_{n,m}$ as \textit{generalized sampling}.  Observe that $\cP_{m} f$ is determined solely by the coefficients $\hat f_1,\ldots \hat f_m$ of $f$.  Hence $F_{n,m}(f)$ is also determined solely by these values.

In what follows it will be useful to note that \R{redconsisteqns} is equivalent to
\be{
\label{reconsisteqns2}
\ip{\tilde{f}_{n,m}}{\cP_{m} \phi_j} = \ip{f }{\cP_{m} \phi_j},\quad j=1,\ldots,n,\qquad \tilde f_{n,m} \in \rT_n,
}
due to the self-adjointness of $\cP_{m}$.  An immediate consequence of this formulation is the following:

\lem{
\label{l:consistGS} Suppose that $\cos (\theta_{n,n}) >0$ and that $\dim(\rS_n) = \dim(\rT_n)$.  Then when $m =n $ the generalized sampling reconstruction $\tilde f_{n,m}$ of $f \in \rH$ defined by \R{redconsisteqns} is precisely the consistent reconstruction $\tilde f_{n,n}$ defined by \R{E:findimconsist}.
}
\prf{
We first claim that $\cP_n$ is a bijection from $\rT_n$ to $\rS_n$.  Suppose that $\cP_n \phi = 0$ for some $\phi \in \rT_n$.  Then $0 = \ip{\cP_n \phi}{\phi} = \sum^{n}_{j=1} | \ip{\phi}{\psi_j} |^2$ and therefore $\phi \in \rS^{\perp}_ {n}$.  Since $\phi \in \rT_n$, and $\rT_n \cap \rS^{\perp}_{n} = \{ 0 \}$ by assumption, we have $\phi = 0$, as required.

By linearity, we now find that the conditions \R{reconsisteqns2} are equivalent to \R{E:findimconsist}.  Since the consistent reconstruction $\tilde f_{n,n}$ satisfying \R{consistconds} exists uniquely (Corollary \ref{c:finconsist}), we obtain the result.
}

We conclude that generalized sampling contains consistent sampling as a special case corresponding to $n=m$, which justifies the use of the same notation for both.  However, as expounded in \cite{BAACHShannon,BAACHAccRecov}, the key to generalized sampling is to allow $m$ to vary independently from $n$.  As we prove, doing so results in a small reconstruction constant.

\subsection{An intuitive argument}
Before providing a complete analysis, let us first give an intuitive explanation as to why this is the case.  To this end, suppose that $n$ is fixed and let $m \rightarrow \infty$.  Equations \R{redconsisteqns} now read
\bes{
\ip{\cP\tilde{f}_{n,\infty}}{\phi}= \ip{\cP f}{\phi},\quad \forall \phi \in \rT_{n},\qquad \tilde f_{n,\infty} \in \rT_n,
}
for some $\tilde{f}_{n,\infty} \in \rT_{n}$, where $\cP$ is given by \R{Pdef}.  In \cite{BAACHAccRecov} it was shown that $\tilde{f}_{n,\infty} = \lim_{m \rightarrow \infty} \tilde{f}_{n,m}$ for fixed $n \in \bbN$, exactly as one would expect.  Hence, we can understand the behaviour of $\tilde{f}_{n,m}$ for large $m$ by first analysing $\tilde{f}_{n,\infty}$.

Since $\cP$ is self-adjoint, $\tilde f_{n,\infty}$ is equivalently defined by
\be{
\label{E:fntilde}
\ip{\tilde f_{n,\infty}}{\Phi} = \ip{f}{\Phi},\quad \forall \Phi \in \cP(\rT_n),\qquad \tilde f_{n,\infty} \in \rT_n.
}
We now have

\thm{
\label{t:intuitive}	
For any $f \in \rH$, there exists a unique $\tilde f_{n,\infty} \in \rT_n$ satisfying \R{E:fntilde}.  Moreover, the mapping $f \mapsto \tilde f_{n,\infty}$ is precisely the oblique projection with range $\rT_n$ and kernel $(\cP(\rT_n))^{\perp}$, and we have the sharp bounds
\be{
\label{E:fnstab}
\| \tilde f_{n,\infty} \| \leq \sec \left ( \theta_{n,\infty} \right ) \| f \|,
}
and
\be{
\label{E:fnerr}
\| f - \cQ_n f \| \leq \| f - \tilde f_{n,\infty} \| \leq \sec \left ( \theta_{n,\infty} \right ) \| f - \cQ_{n} f \|,
}
where $ \theta_{n,\infty} $ is the angle between $\rT_n$ and $\cP(\rT_n)$.
}
\prf{
We first claim that  $\cos(\theta_{n,\infty}) > 0$, so that the oblique projection $\cW$ with range $\rT_n$ and kernel $(\cP(\rT_n))^{\perp}$ is well-defined as a mapping of $\rH_ 0 = \rT_n \oplus (\cP(\rT_n))^{\perp}$.  Suppose not.  Since $\rT_n$ is finite dimensional, there exists a $\phi \in \rT_n$, $\phi \neq 0$, satisfying $\cQ_{\cP(\rT_n)} \phi = 0$.  Thus
\bes{
0 = \ip{\cQ_{\cP(\rT_n)} \phi}{\cP \phi'} = \ip{\phi}{\cP \phi'},\quad \forall \phi' \in \rT_n,
}
and, in particular, $\ip{\phi}{\cP \phi} = 0$.  Thus $\phi = 0$ by Lemma \ref{l:Pbdd} -- a contradiction.  Hence $\cW$ is well-defined.  Note that $\cW f$ satisfies the equations \R{E:fntilde}.  Arguing in the standard way, we can show that solutions to \R{E:fntilde} are unique.  Hence $\tilde{f}_{n,\infty} = \cW f$, as required.

It remains to show that $\rH_0 = \rH$.  The result follows immediately from Lemma \ref{l:findimoblique} provided $\dim(\cP(\rT_n)) = \dim(\rT_n)$.  However, if not, then there is a nonzero $\phi \in \rT_n$ with $\cP \phi = 0$, which also contradicts Lemma \ref{l:Pbdd}.
}

Note that this theorem improves on \cite[Thm.\ 2.1]{BAACHAccRecov} by giving sharp bounds.  We can also estimate the reconstruction constant of the mapping $F_{n,\infty} : f \mapsto \tilde{f}_{n,\infty}$:
\cor{
\label{c:intuitive}
Let $F_{n,\infty}$ be the mapping $f \mapsto \tilde{f}_{n,\infty}$, where $\tilde{f}_{n,\infty}$ is defined by \R{E:fntilde}.  Then 
\be{
\label{stat3inf}
1 \leq  \mu(F_{n,\infty}) \leq \frac{\sqrt{c_2}}{\sqrt{c_1} \cos \left ( \theta_{\rT \rS} \right )},\qquad \frac{1}{\sqrt{c_2}} \leq \kappa(F_{n,\infty}) \leq \frac{1}{\sqrt{c_1} \cos \left ( \theta_{\rT \rS} \right )}. 
}
and
\be{
\label{stat4inf}
\max \left \{ 1 , \frac{1}{\sqrt{c_2}} \right \} \leq C(F_{n,\infty}) \leq \frac{ \max \left \{ 1 , \sqrt{c_2} \right \} }{\sqrt{c_1} \cos \left ( \theta_{\rT \rS} \right )}.
}
}
This corollary (we present the proof in the next section) confirms that the reconstruction scheme $\{ F_{n,\infty} \}_{n \in \bbN}$ is stable and quasi-optimal in the sense of Definition \ref{d:reconschemestab}, with constant $C^* \leq \frac{\max\{1, \sqrt{c_2} \}}{\sqrt{c_1} \cos (\theta_{\rT \rS})}$.  Hence, unlike the consistent sampling scheme $\{ F_{n,n} \}_{n \in \bbN}$, where the reconstruction constants $C(F_{n,n})$ can quite easily be exponentially large in $n$ (see \S \ref{ss:consistfail}), this approach guarantees good approximation and robustness with respect to noise.

However, one cannot actually compute $\tilde{f}_{n,\infty}$, since it involves the infinite-rank operator $\cP$.  Nonetheless, since the generalized sampling reconstruction $\tilde{f}_{n,m} \approx \tilde{f}_{n,\infty}$ for large $m$, we may expect the good properties of $F_{n,\infty}$, i.e.\ stability and quasi-optimality, to be inherited whenever $m$ is sufficiently large.  In the next section we prove this to be the case.

Before doing so, however, let us relate $\tilde f_{n,\infty}$ to the the discussion in \S \ref{ss:optheory}.  Recall first that we wish to solve $U \alpha = \hat{f}$.  Since $\alpha$ satisfies these equations it also obeys the normal equations
\be{
\label{E:normal}
U^* U \alpha = U^* \hat{f}.
}
Now write $\tilde f_{n,\infty} = \sum^{n}_{j=1} \alpha^{[n,\infty]}_j \phi_j$.  It is easily shown that $\alpha^{[n,\infty]} = \{ \alpha^{[n,\infty]}_1,\ldots,\alpha^{[n,\infty]}_n \}$ is defined by
\bes{
P_{n} U^* U P_n \alpha = P_n U^* \hat{f},
}
where $\hat{f} = \{ \hat{f}_1,\hat{f}_2,\ldots \}$.  Thus, $\alpha^{[n]}$ is precisely the result of the finite section method applied to the normal equations \R{E:normal}.  Since the operator $U^* U$ is self-adjoint and positive, its finite sections must possess properties (i)--(iii), and hence we are guaranteed a good reconstruction.

\subsection{Analysis of generalized sampling}
The analysis of $\tilde f_{n,m}$ is similar to that of $\tilde f_{n,\infty}$.  Whilst such an analysis was originally given in \cite{BAACHShannon,BAACHAccRecov}, the estimates derived were not sharp.  Our main result in this section is to present new, sharp bounds.  We first require the following lemma:
\lem{
\label{l:theta_conv}
Let $\theta_{n,m}$ and $\theta_{n,\infty}$ be the angles between $\rT_n$ and the subspaces $\cP_m(\rT_n)$ and $\cP(\rT_n)$ respectively.  Then, for fixed $n$, $\theta_{n,m} \rightarrow \theta_{n,\infty}$ as $m \rightarrow \infty$.  In particular,
\bes{
1 \leq \lim_{m \rightarrow \infty} \sec \left( \theta_{n,m} \right ) \leq \frac{\sqrt{c_2}}{\sqrt{c_1} \cos \left ( \theta_{\rT \rS} \right )}.
}
}
\prf{
From the definition \R{subangle}, we have
\bes{
\cos \left ( \theta_{n,m} \right ) = \inf_{\substack{\phi \in \rT_n \\ \| \phi \| = 1}} \sup_{\substack{\phi' \in \rT_n \\ \cP_m \phi' \neq 0}} \frac{\ip{\phi}{\cP_m \phi'}}{\| \cP_m \phi' \|}.
}
Recall first that $\cP_m \rightarrow \cP$ strongly on $\rH$.  Since $\rT_n$ is finite-dimensional, this implies uniform convergence of $\cP_m \rightarrow \cP$ on $\rT_n$.  In other words, if $e_m = \| \cP |_{\rT_n}  - \cP_m |_{\rT_n} \|$ then $e_m \rightarrow 0$ as $m \rightarrow \infty$.  In particular, for sufficiently large $m$, $\cP_m \phi' \neq 0$ if and only if $\cP \phi' \neq 0$.  Thus, for large $m$,
\be{
\label{skip0}
\cos \left ( \theta_{n,m} \right ) = \inf_{\substack{\phi \in \rT_n \\ \| \phi \| = 1}} \sup_{\substack{\phi' \in \rT_n \\  \phi' \neq 0}} \frac{\ip{\phi}{\cP_m \phi'}}{\| \cP_m \phi' \|} = \inf_{\substack{\phi \in \rT_n \\ \| \phi \| = 1}} \sup_{\substack{\phi' \in \rT_n \\ \cP \phi' \neq 0}} \frac{\ip{\phi}{\cP_m \phi'}}{\| \cP_m \phi' \|}.
}
Now
\ea{
\label{skip1}
\frac{\ip{\phi}{\cP_m \phi'}}{\| \cP_m \phi' \|} &= \left ( \frac{ \| \cP \phi' \|}{\| \cP_m \phi' \|} \right ) \left ( \frac{\ip{\phi}{\cP \phi'} }{ \| \cP \phi' \|} - \frac{\ip{\phi}{(\cP - \cP_m) \phi'}}{ \| \cP \phi' \|}\right ).
}
Note that
\bes{
\big | \| \cP_m \phi' \| - \| \cP \phi' \| \big | \leq \| (\cP - \cP_m) \phi' \| \leq e_m \| \phi' \|.
}
Moreover,  
\bes{
\| \cP \phi' \| = \sup_{\substack{g \in \rH \\ \| g \|=1}} \ip{\cP \phi'}{g} \geq \frac{\ip{\cP \phi'}{\phi'}}{\| \phi'\|} \geq c_1 \cos^2 \left ( \theta_{\rT \rS} \right ) \| \phi' \|.
}
Thus,
\be{
\big | \| \cP_m \phi' \| - \| \cP \phi' \| \big | \leq \frac{e_m}{c_1} \sec^2 ( \theta_{\rT \rS} ) \| \cP \phi' \|.
}
Combining this with \R{skip1}, we now obtain
\bes{
\frac{\ip{\phi}{\cP_m \phi'}}{\| \cP_m \phi' \|} \geq \frac{1}{1+ \frac{e_m}{c_1} \sec^2 ( \theta_{\rT \rS} )} \left ( \frac{\ip{\phi}{\cP \phi'} }{ \| \cP \phi' \|} - \| \phi \| \frac{e_m}{c_1} \sec^2 ( \theta_{\rT \rS} ) \right ),
}
and
\bes{
\frac{\ip{\phi}{\cP_m \phi'}}{\| \cP_m \phi' \|} \leq \frac{1}{1- \frac{e_m}{c_1} \sec^2 ( \theta_{\rT \rS} )} \left ( \frac{\ip{\phi}{\cP \phi'} }{ \| \cP \phi' \|} + \| \phi \| \frac{e_m}{c_1} \sec^2 ( \theta_{\rT \rS} ) \right ),
}
Hence \R{skip0} now gives
\eas{
\cos ( \theta_{n,m} ) & \geq \frac{1}{1+ \frac{e_m}{c_1} \sec^2 ( \theta_{\rT \rS} )}  \left ( \cos ( \theta_{n,\infty}  ) - \frac{e_m}{c_1} \sec^2 ( \theta_{\rT \rS} ) \right ) \\
\cos ( \theta_{n,m} ) & \leq  \frac{1}{1- \frac{e_m}{c_1} \sec^2 ( \theta_{\rT \rS} )}  \left ( \cos ( \theta_{n,\infty}  ) + \frac{e_m}{c_1} \sec^2 ( \theta_{\rT \rS} ) \right ),
}
and the result now follows from the fact that $e_m \rightarrow 0$ as $m \rightarrow \infty$.
}

We now have:

\thm{
\label{t:reconthm}
For each $n \in \bbN$ and any $f \in \rH$, there exists an $m_{0}$, independent of $f$, such that the reconstruction $\tilde{f}_{n,m}$ defined by \R{redconsisteqns} exists and is unique for all $m \geq m_{0}$.  In particular, $m_0$ is the least $m$ such that $\cos (\theta_{n,m})>0$.  Moreover, the mapping $f \mapsto \tilde f_{n,m}$ is precisely the oblique projection onto $\rT_n$ along $[\cP_{m}(\rT_n)]^{\perp}$, and we have the sharp bounds
\be{
\label{E:fnmstab}
\| \tilde f_{n,m} \| \leq \sec \left ( \theta_{n,m} \right ) \| f \|,
}
and
\be{
\label{E:fnmerr}
\| f - \cQ_n f \| \leq \| f - \tilde f_{n,m} \| \leq \sec \left ( \theta_{n,m} \right ) \| f - \cQ_{n} f \|.
}
}
\prf{
The existence of an $m_0$ such that $\cos(\theta_{n,m}) >0$ for all $m \geq m_0$ follows from Lemma \ref{l:theta_conv}.  Thus, when $m \geq m_0$ the oblique projection $\cW$ with range  $\rT_n$ and kernel $(\cP_{m}(\rT_n))^{\perp}$ is well-defined over $\rH_ 0 := \rT_n \oplus (\cP_{m}(\rT_n))^{\perp}$ and satisfies \R{redconsisteqns}.  Moreover, under this condition, solutions of \R{redconsisteqns} are unique, and thus $\tilde{f}_{n,m} = \cW f$ whenever $f \in \rH_0$.  An application of Corollary \ref{c:consist} now gives \R{E:fnmstab} and \R{E:fnmerr}.  To complete the proof we need only show that $\rH_ 0 = \rH$.  This follows immediately from Corollary \ref{l:findimoblique}, provided $\dim ( \cP_m(\rT_n ) ) = \dim (\rT_n)$.  However, if not then there exists a nonzero $\phi \in \rT_n \cap (\cP_m(\rT_n))^{\perp}$, which contradicts the fact that $\cos (\theta_{n,m} ) > 0$.}

Note that this theorem improves the bounds of \cite[Thm.\ 2.4]{BAACHAccRecov}, and gives the exact value $\mu(F_{n,m} ) = \sec (\theta_{n,m})$ for the quasi-optimality constant  of generalized sampling.  Having done this, we next determine the condition number $\kappa(F_{n,m})$, and as a result the reconstruction constant $C(F_{n,m})$.  For this, we introduce the following quantity:
\be{
\label{Dnm_def}
D_{n,m} = \left ( \inf_{\substack{\phi \in \rT_n \\ \| \phi \|=1}} \ip{\cP_m \phi}{\phi} \right )^{-\frac12},\quad n,m\in \bbN,
}
(this is similar to the quantity $C_{n,m}$ of \cite[Eqn.\ (2.12)]{BAACHAccRecov}).  Note that $D_{n,m}$ need not be defined for all $n,m\in \bbN$.  However, we will show subsequently that this is the case provided $m$ is sufficiently large (for a given $n$).  We shall also let
\bes{
D_{n,\infty} = \left ( \inf_{\substack{\phi \in \rT_n \\ \| \phi \|=1}} \ip{\cP \phi}{\phi} \right )^{-\frac12},\quad n\in \bbN.
}
We now have the following lemma:
\lem{
\label{l:Dnmconv}
For fixed $n \in \bbN$, $D_{n,m} \rightarrow D_{n,\infty}$ as $m \rightarrow \infty$.  In particular, 
\bes{
 \frac{1}{\sqrt{c_2}} \leq \lim_{m \rightarrow \infty} D_{n,m} \leq \frac{1}{\sqrt{c_1} \cos \left ( \theta_{\rT \rS} \right )}.
}
}
\prf{
The first result follows from strong convergence of the operators $\cP_m \rightarrow \cP$ on $\rH$ and the fact that $\rT_n$ is finite-dimensional.  The second result is due to Lemma \ref{l:Pbdd}.
}

With this to hand, our main result is as follows:

\cor{
\label{c:GSreconconst}
Let $n \in \bbN$ and suppose that $m \geq m_0$, where $m_0$ is as in Theorem \ref{t:reconthm}.  Let $F_{n,m}$ be the generalized sampling reconstruction $f \mapsto \tilde{f}_{n,m}$, where $\tilde{f}_{n,m}$ is defined by \R{redconsisteqns}.  Then 
\be{
\label{stat1}
\mu(F_{n,m}) = \sec \left (\theta_{n,m} \right) ,\quad \kappa(F_{n,m}) = D_{n,m},
}
and the reconstruction constant $C(F_{n,m})$ satisfies
\be{
\label{stat2}
D_{n,m} \leq C(F_{n,m}) \leq \max \left \{ 1 , \sqrt{c_2} \right \} D_{n,m},
}
whenever $D_{n,m}$ is defined.  In particular, for fixed $n$,
\be{
\label{stat3}
1 \leq \lim_{m \rightarrow \infty} \mu(F_{n,m}) \leq \frac{\sqrt{c_2}}{\sqrt{c_1} \cos \left ( \theta_{\rT \rS} \right )},\qquad \frac{1}{\sqrt{c_2}} \leq \lim_{m \rightarrow \infty} \kappa(F_{n,m}) \leq \frac{1}{\sqrt{c_1} \cos \left ( \theta_{\rT \rS} \right )},
}
and
\be{
\label{stat4}
\max \left \{ 1 , \frac{1}{\sqrt{c_2}} \right \} \leq \lim_{m \rightarrow \infty} C(F_{n,m}) \leq \frac{ \max \left \{ 1 , \sqrt{c_2} \right \} }{\sqrt{c_1} \cos \left ( \theta_{\rT \rS} \right )}.
}
}
\prf{
We claim that
\be{
\label{Pclaim}
\sec \left ( \theta_{n,m} \right ) \leq \sqrt{c_2} D_{n,m}
}
Note first that $D_{n,m} < \infty$ implies that $\cP_m |_{\rT_n} : \rT_n \rightarrow \cP_m(\rT_n)$ is invertible.  Hence, by \R{subangle},
\be{
\label{cosineq}
\cos \left ( \theta_{n,m} \right ) = \inf_{\substack{\phi \in \rT_n \\ \phi \neq 0 }} \sup_{\substack{\phi' \in \rT_n \\ \phi' \neq 0 }} \frac{\ip{\phi}{\cP_m \phi'}}{\| \phi \| \| \cP_m \phi ' \|} \geq \inf_{\substack{\phi \in \rT_n \\ \phi \neq 0 }} \frac{\ip{\phi}{\cP_m \phi}}{\| \phi \| \| \cP_m \phi \|}.
}
Now consider $\| \cP_m \phi \|$.  Since $\cP_m \phi \in \rS_m$, basic properties of $\cP_m$ give that
\ea{
\label{PmNormUB}
\| \cP_m \phi \| = \sup_{\substack{\psi \in \rS_m \\ \| \psi \| = 1}} \ip{\psi}{\cP_m \phi}  \leq \sqrt{\ip{\phi}{\cP_m \phi}} \sup_{\substack{\psi \in \rS_m \\ \| \psi \| = 1}} \sqrt{\ip{\psi}{\cP_m \psi}}  \leq \sqrt{\ip{\phi}{\cP_m \phi}} \sup_{\substack{\psi \in \rS \\ \| \psi \| = 1}} \sqrt{\ip{\psi}{\cP \psi}}, 
}
and therefore 
\bes{
\| \cP_m \phi \| \leq \sqrt{c_2} \sqrt{\ip{\phi}{\cP_m \phi}}.
}
Applying this to \R{cosineq} now gives  \R{Pclaim}.

Note that \R{stat2} and \R{stat4} now follow immediately from \R{stat1},  \R{stat3}, \R{Pclaim} and the definition $C(F_{n,m}) = \max \{ \mu(F_{n,m}) , \kappa(F_{n,m}) \}$.  Moreover, \R{stat3} follows from \R{stat1} and Lemmas \ref{l:theta_conv} and \ref{l:Dnmconv}.  Hence we only need to prove \R{stat1}.  The first part is due to Theorem \ref{t:reconthm}.  Therefore it remains to show that $\kappa(F_{n,m}) = D_{n,m}$.  Since $F_{n,m}$ is linear and, due to \R{E:fnmerr}, perfect on $\rT_n$, we have
\bes{
\kappa(F_{n,m}) = \sup_{\substack{f \in \rH \\ \hat{f} \neq 0}} \left \{ \frac{\| F_{n,m}(f) \|}{\| \hat{f} \|_{\ell^2}} \right \} \geq \sup_{\substack{\phi \in \rT_n \\ \phi \neq 0}} \left \{ \frac{\| \phi \|}{\| \hat{\phi} \|_{\ell^2} } \right \}.
}
Since $\| \hat{\phi} \|^2_{\ell^2} = \ip{\cP_m \phi}{\phi}$, this now gives
\bes{
\kappa(F_{n,m}) \geq \left ( \inf_{\substack{\phi \in \rT_n \\ \| \phi \| = 1}} \ip{\cP_m \phi}{\phi} \right )^{-\frac12} = D_{n,m}.
}
We now wish to derive the upper bound.  Let $F_{n,m}(f) = \tilde{f}_{n,m}$.  Since $\tilde{f}_{n,m} \in \rT_n$, \R{redconsisteqns} gives that
\bes{
\ip{\cP_m \tilde{f}_{n,m}}{\tilde{f}_{n,m}} = \ip{\cP_m f}{\tilde{f}_{n,m}} \leq \sqrt{\ip{\cP_m f}{f}} \sqrt{ \ip{\cP_m \tilde{f}_{n,m}}{\tilde{f}_{n,m}}}.
}
Thus $\| \hat{f} \|_{\ell^2} \geq \sqrt{ \ip{\cP \tilde{f}_{n,m}}{\tilde{f}_{n,m}}}$.  Since $f \mapsto \tilde{f}_{n,m}$ is a surjection onto $\rT_n$, we therefore deduce that
\bes{
\kappa(F_{n,m} ) \leq \sup_{\substack{f \in \rH \\ \hat{f} \neq 0}} \left \{ \frac{\| \tilde{f}_{n,m} \|}{\sqrt{ \ip{\cP_m \tilde{f}_{n,m}}{\tilde{f}_{n,m}}}} \right \} = \sup_{\substack{\phi \in \rT \\ \phi \neq 0 }} \left \{ \frac{\| \phi \|}{\sqrt{\ip{\cP_m \phi}{\phi}}} \right \} = D_{n,m},
}
as required.  This completes the proof.
}

We are now in a position to establish the two results not proved previously: 
\prf{[Proof of Corollary \ref{c:intuitive}]
By replacing $\cP_m$ by $\cP$ in the proof of Corollary \ref{c:GSreconconst} we find that $\sec \left ( \theta_{n,\infty} \right ) \leq \sqrt{c_2} D_{n,\infty}$ and $\kappa(F_{n,\infty}) = D_{n,\infty}$.  The result now follows from this and Lemmas \ref{l:theta_conv} and \ref{l:Dnmconv}. 
}

\prf{[Proof of Corollary \ref{c:reconconst}]
Under the assumption $\dim(\rS_n) = \dim(\rT_n)$ the consistent reconstruction \R{E:findimconsist} coincides with the generalized sampling reconstruction \R{redconsisteqns} (Lemma \ref{l:consistGS}).  The result now follows immediately from \R{stat1} and \R{Pclaim}.
}

Corollary \ref{c:GSreconconst} confirms the advantage of generalized sampling.  Given $n \in \bbN$, one can always take $m$ sufficiently large to guarantee a stable, quasi-optimal reconstruction with reconstruction constant asymptotically bounded by $\frac{ \max \left \{ 1 , \sqrt{c_2} \right \} }{\sqrt{c_1} \cos \left ( \theta_{\rT \rS} \right )}$.  

The key issue remaining is to determine how large $m$ must be taken to ensure such properties.  This will be discussed in the next section.  First, however, let us connect generalized sampling to the discussion of \S \ref{ss:optheory}.  Observe that if $\tilde f_{n,m} = \sum^{n}_{j=1} \alpha^{[n,m]}_{j} \phi_j$, then the vector $\alpha^{[n,m]} = \{ \alpha^{[n,m]}_{1},\ldots,\alpha^{[n,m]}_{n} \} \in \bbC^n$ is the unique solution to
\bes{
(U^{[n,m]})^{*} U^{[n,m]} \alpha^{[n,m]} = (U^{[n,m]})^* \hat{f}^{[m]},
}
where $U^{[n,m]} \in \bbC^{m \times n}$ is precisely $P_{m} U P_{n}$ and $\hat{f}^{[m]} = P_m \hat{f}$.  The matrix $U^{[n,m]}$ is the leading $m \times n$ submatrix of $U$, and is sometimes referred to as an \textit{uneven section} of $U$.  Uneven sections have recently gained prominence as effective alternatives to the finite section method for discretizing non-self adjoint operators \cite{strohmer,Lindner2008}.  In particular, in \cite{hansen2011} they were employed to solve the long-standing computational spectral problem.  Their success is due to the observation that, under a number of assumptions (which are always guaranteed for the problem we consider in this paper),  we have
\bes{
(U^{[n,m]})^{*} U^{[n,m]} = P_n U^* P_m U P_n \rightarrow P_n U^* U P_n,\quad m \rightarrow \infty,
}
where $P_n U^* U P_n$ is the $n \times n$ finite section of the self-adjoint matrix $U^* U$.  This guarantees properties (i)--(iii) for $U^{[n,m]}$, whenever $m$ is sufficiently large in comparison to $n$.

Finite (and uneven) sections have been extensively studied \cite{bottcher1996,hansen2008,lindner2006}, and there exists a well-developed and intricate theory of their properties involving $C^{*}$-algebras \cite{hagen2001c}.  However, these general results say little about the rate of convergence, nor do they provide explicit constants.  Yet, as illustrated in Theorem \ref{t:reconthm}, the operator $U$ in this case is so structured that its uneven sections admit both explicit constants and estimates for the rate of convergence.  Moreover, of great importance for computations, such constants can also be numerically computed, as we discuss in \S \ref{ss:stabsamp}.

\subsection{The condition number and quasi-optimality constant}\label{ss:theta_D_relation}
As shown, the condition number $\kappa(F_{n,m})$ coincides with $D_{n,m}$, and the quasi-optimality constant $\mu(F_{n,m})=\sec \left ( \theta_{n,m} \right)$, where $\theta_{n,m}$ is the angle between $\rT_n$ and $\cP_m(\rT_n)$.  In addition, since
\be{
\label{sec_theta_bound}
\sec \left ( \theta_{n,m} \right) \leq \sqrt{c_2} D_{n,m},
}
one can control the behaviour of both quantities, and therefore also $C(F_{n,m})$, by controlling $D_{n,m}$.  The advantage of this, as we discuss in \S \ref{ss:stabsamp}, is that it is typically easier to compute $D_{n,m}$ than it is $\theta_{n,m}$.

However, it is in general possible for $\sec \left ( \theta_{n,m} \right )$ to be somewhat smaller than $D_{n,m}$.  Thus, in numerical examples, one may see a better approximation than the bound \R{sec_theta_bound} suggests.  For this reason, we now present a result that addresses the relationship between $\sec \left ( \theta_{n,m} \right)$ and $D_{n,m}$:

\lem{
\label{l:sec_D_relation}
Let $G^{[m]} = \{ \ip{\psi_j}{\psi_k} \}^{m}_{j,k=1} \in \bbC^{m \times m}$ be the Gram matrix of the first $m$ sampling vectors $\{\psi_1,\ldots,\psi_m\}$.  Then
\bes{
\sqrt{c_{1,m}} D_{n,m} \leq \sec \left ( \theta_{n,m} \right) \leq \sqrt{c_{2,m}} D_{n,m},
}
where $c_{1,m}$ and $c_{2,m}$ are the frame bounds for the frame sequence $\{ \psi_1,\ldots,\psi_m \}$.  In particular, when $\{ \psi_j \}_{j \in \bbN}$ is a Riesz basis, we have
\bes{
d_1 D_{n,m} \leq \sec \left ( \theta_{n,m} \right) \leq d_2 D_{n,m},
}
where $d_1,d_2>0$ are the Riesz basis constants for $\{ \psi_j \}_{j \in \bbN}$, and when $\{ \psi_j \}_{j \in \bbN}$ is an orthonormal basis it holds that
\bes{
\sec \left ( \theta_{n,m} \right) = D_{n,m}.
}
}  
Note that, by Riesz basis constants, we mean constants $d_1,d_2>0$ such that
\bes{
d_1 \| \beta \|_{\ell^2} \leq \nm{\sum_{j \in \bbN} \beta_j \psi_j } \leq d_2 \| \beta \|_{\ell^2},\quad \forall \beta = \{ \beta_j \}_{j \in \bbN} \in \ell^2(\bbN).
}

\prf{[Proof of Lemma \ref{l:sec_D_relation}]
By definition
\bes{
\cos \left ( \theta_{n,m} \right ) = \inf_{\substack{\phi \in \rT_n \\ \| \phi \|=1}} \| \cQ_{\cP_m(\rT_n)} \phi \|.
}
Consider $\| \cQ_{\cP_m(\rT_n)} \phi \|$.  Since $\cP_m(\rT_n) \subseteq \rS_m$, we have
\bes{
\| \cQ_{\cP_m(\rT_n)} \phi \| = \sup_{\substack{\psi \in \cP_m(\rT_n) \\ \| \psi \|=1}} \ip{\phi}{\psi} \leq \sup_{\substack{\psi \in \rS_m \\ \| \psi \|=1}} \ip{\phi}{\psi}.
}
Recall that the operator $\cP_m$ is invertible on $\rS_m$.  Hence
\be{
\label{jump1}
\| \cQ_{\cP_m(\rT_n)} \phi \| \leq \sup_{\substack{\psi \in \rS_m \\ \psi \neq 0}} \frac{ \ip{\phi}{\cP_m\psi} }{\| \cP_m \psi \|}\leq \sqrt{\ip{\cP_m \phi}{\phi}} \sup_{\substack{\psi \in \rS_m \\ \psi \neq 0}} \frac{\sqrt{\ip{\cP_m \psi}{\psi}}}{\| \cP_m \psi \|}.
}
Consider the latter term.  The operator $\cP_m : \rS_m \rightarrow \rS_m$ is invertible, self-adjoint and positive-definite.  Hence, it has a unique square root $(\cP_m)^{\frac12}$ with these properties \cite[Lem. 2.4.4]{christensen2003introduction}.  Thus
\bes{
\sup_{\substack{\psi \in \rS_m \\ \psi \neq 0}} \frac{\sqrt{\ip{\cP_m \psi}{\psi}}}{\| \cP_m \psi \|} = \sup_{\substack{\psi \in \rS_m \\ \psi \neq 0}} \frac{\| (\cP_m)^{\frac12} \psi \| }{\| \cP_m \psi \|} = \sup_{\substack{\psi \in \rS_m \\ \psi \neq 0}} \frac{\| \psi \| }{\| (\cP_m)^{\frac12}  \psi \|} = \sup_{\substack{\psi \in \rS_m \\ \psi \neq 0}} \frac{\| \psi \| }{\sqrt{\ip{\cP_m \psi}{\psi}}} = \frac{1}{\sqrt{c_{1,m}}}.
}
Combining this with \R{jump1} now gives $\sec ( \theta_{n,m}) \geq \sqrt{c_{1,m}} D_{n,m}$ as required.

For the upper bound, we first notice that
\bes{
\| \cQ_{\cP_m(\rT_n)} \phi \| \geq \frac{\ip{\phi}{\cP_m \phi}}{\| \cP_m \phi \|}.
}
Moreover, arguing as in \R{PmNormUB} one finds that $\| \cP_m \phi \| \leq \sqrt{c_{2,m}} \sqrt{\ip{\cP_m \phi}{\phi}}$.  Combining this with the previous expression and using the definition of $\cos(\theta_{n,m})$ now gives the first result.

For the second part of the proof, we first recall that the Riesz basis bounds $d_1,d_2$ for $\{ \psi_j \}_{j \in \bbN}$ are lower and upper bounds for the Riesz basis bounds $d_{1,m},d_{2,m}$  for the finite subset $\{ \psi_1,\ldots,\psi_m \}$.  Moreover, by \cite[Thm. 5.2.1]{christensen2003introduction}, the frame bounds $c_{1,m},c_{2,m}$ for the Riesz basis $\{ \psi_1,\ldots,\psi_m \}$ are identical to the Riesz basis bounds $d_{1,m},d_{2,m}$.  This gives the second result.  Finally, when $\{ \psi_j \}_{j \in \bbN}$ is orthonormal we have $d_1=d_2 =1$, and thus we obtain the final result.
}
This lemma demonstrates that the difference in magnitudes between $\sec \left ( \theta_{n,m} \right )$ and $D_{n,m}$ is determined by $\sqrt{c_{1,m}}$ and $\sqrt{c_{2,m}}$.  Note that $c_{2,m} \leq c_2$, where $c_2$ is the frame bound for the infinite frame $\{ \psi_j \}_{j \in \bbN}$.  However, $c_{1,m}$ can exhibit wild behaviour: it is possible to construct simple frames for which $c_{1,m}$ is exponentially small in $m$, even though $c_1$ is moderate in magnitude \cite{Christensen1993Projection}.  On the other hand, if $\{ \psi_j \}_{j \in \bbN}$ is a Riesz or orthonormal basis, we find that $D_{n,m}$ and $\sec (\theta_{n,m} )$ are, up to a possible factor proportional to the Riesz basis constants $d_1$ and $d_2$, the same.

\subsection{Computing the generalized sampling reconstruction}\label{ss:GScompute}
Recall that the generalized sampling reconstruction $\tilde{f}_{n,m}$ depends only on $\rT_n$, and not on the system of functions used to span $\rT_n$.  Let $\{ \phi_j \}^{d_n}_{j=1}$ be a spanning set for $\rT_n$, where $d_n \geq \dim(\rT_n)$, and write
\bes{
\tilde{f}_{n,m} = \sum^{d_n}_{j=1} \alpha^{[n,m]}_{j} \phi_j.
}
The vector $\alpha^{[n,m]} = \{ \alpha^{[n,m]}_{j} \}^{d_n}_{j=1}$ is the least squares solution to the overdetermined linear system
\bes{
U^{[n,m]} \alpha^{[n,m]} = \hat{f}^{[m]}, 
}
where $U^{[n,m]} \in \bbC^{m \times d_n}$ has $(j,k)^{\rth}$ entry $\ip{\phi_k}{\psi_j}$.  Thus, computing $\tilde{f}_{n,m}$ is equivalent to solving a least squares problem.  From a numerical perspective, it is important to understand the condition number $\kappa(U^{[n,m]}) = \| U^{[n,m]} \| \| (U_{n,m})^{\dag} \|$ of the matrix $U^{[n,m]}$, where $\dag$ denotes the pseudoinverse.  The following lemma is similar to \cite[Lem.\ 2.11]{BAACHAccRecov} (for this reason we forgo the proof):
\lem{
Let $\{ \phi_j \}^{d_n}_{j=1}$ be a spanning set for $\rT_n$, and write $G^{[n]} \in \bbC^{d_n \times d_n}$ for its Gram matrix.  Then the condition number of the matrix $U^{[n,m]}$ satisfies
\bes{
\frac{1}{\sqrt{c_2} D_{n,m}} \sqrt{\kappa \left(G^{[n]} \right )} \leq \kappa(U^{[n,m]}) \leq \sqrt{c_2} D_{n,m} \sqrt{\kappa \left(G^{[n]} \right )}.
}
}
This lemma shows that the condition number of the matrix $U^{[n,m]}$ is no worse than that of the Gram matrix $G^{[n]}$ whenever $m$ is chosen sufficiently large to ensure boundedness of $D_{n,m}$.  In particular, if the vectors $\{ \phi_1,\ldots,\phi_n \}$ are a Riesz or orthonormal basis, then $\kappa(G^{[n]}) = \ord{1}$ and hence the condition number of $U^{[n,m]}$ is completely determined by the magnitude of $D_{n,m}$.  In this case, not only is the reconstruction $\tilde{f}_{n,m}$ numerically stable, but so is the computation of its coefficients $\alpha^{[n,m]}$.  For further details on the computation of $\tilde{f}_{n,m}$, see \cite{BAACHAccRecov}.

\section{The stable sampling rate}\label{ss:stabsamp}
The key ingredient of generalized sampling is that the parameter $m$ must be sufficiently large in comparison to $n$.  The notion of how large was first quantified in \cite{BAACHShannon,BAACHAccRecov}.  In this section we improve on this by using the sharp bounds of the previous section.  We define:

\defn{
\label{d:SSR}
The stable sampling rate is defined by
\bes{
\Theta(n;\theta) = \min \left \{ m \in \bbN : C(F_{n,m}) \leq \theta \right \},\quad n \in \bbN,\  \theta \in \left (\frac{\max\{1,\sqrt{c_2} \} }{\sqrt{c_1} \cos (\theta_{\rT \rS})},\infty \right ).
}
}
The stable sampling rate measures how large $m$ must be (for a given $n$) to ensure guaranteed, stable and quasi-optimal recovery.  Indeed,  choosing $m \geq \Theta(n;\theta)$, we find that $C(F_{n,m}) \leq  \theta$ and therefore the reconstruction $\tilde{f}_{n,m}$ is numerically stable and quasi-optimal, up to the magnitude of $\theta$.  In other words, given $n \in \bbN$ and some desired $\theta$, the stable sampling rate determines precisely how many samples are required to guarantee \textit{a priori} a reconstruction constant of magnitude at most $\theta$.  Note that a similar quantity was introduced previously in \cite{BAACHAccRecov}.  However, this was based on estimates for the condition number and quasi-optimality constants which were not sharp.  The stable sampling rate defined above improves on this quantity in that the condition $m \geq \Theta(n;\theta)$ is both sufficient and necessary to ensure stable, quasi-optimal reconstruction: if one were to sample at a rate below the $\Theta(n;\theta)$ then instability and worse convergence of the reconstruction is guaranteed.

One can also ask the reverse question: namely, given a number of samples $m$ and a parameter $\theta$, how large can $n$ be taken?  We refer to the quantity
\be{
\label{SRR}
\Psi(m;\theta) = \max \{ n \in \bbN : C(F_{n,m}) \leq \theta \},\quad m \in \bbN,\  \theta \in \left (\frac{\max\{1,\sqrt{c_2} \} }{\sqrt{c_1} \cos (\theta_{\rT \rS})},\infty \right ),
}
as the \textit{stable reconstruction rate}.

Recall that Remark \ref{r:relax} permits sequences of reconstruction schemes with mildly growing reconstruction constants.  One can also readily define the stable sampling and reconstruction rates to reflect this.  For a positive and increasing sequence $\theta = \{ \theta_n \}_{n \in \bbN}$ with $\inf_{n \in \bbN} \theta_n > \frac{1}{\sqrt{c_1} \cos \left ( \theta_{\rT \rS} \right )}$, we define 
\bes{
\Theta(n ; \theta) = \min \left \{ m \in \bbN : C(F_{n,m}) \leq \theta_n \right \},\quad n \in \bbN,
}
and
\bes{
\Psi(m;\theta) = \max \left \{ n \in \bbN : C(F_{n,m}) \leq \theta_n \right \},\quad m \in \bbN.
}
Once more, one has the interpretation that sampling at the rate $m \geq \Theta(n;\theta)$ ensures stability and quasi-optimality up to the growth of $\theta_n$.

A key property of the stable sampling and reconstruction rates is that they can be computed:

\lem{
Let $D_{n,m}$ and $\theta_{n,m}$ be as in \R{Dnm_def} and Lemma \ref{l:theta_conv} respectively, and suppose that $\{ \phi_j\}^{k_n}_{j=1}$ is a spanning set for $\rT_n$.  Then the quantities $1/D^2_{n,m}$ and $\cos^2 ( \theta_{n,m} )$  are the minimal generalized eigenvalue of the matrix pencils $\left \{ (U^{[n,m]})^* U^{[n,m]} , A^{[n]} \right \}$ and $\{ B^{[n,m]} , A^{[n]} \}$ respectively, where $A^{[n]}$ is the Gram matrix for $\{ \phi_j \}^{k_n}_{j=1}$, $U^{[n,m]}$ is as in \S \ref{ss:GScompute}, $B^{[n,m]}$ is given by
\bes{
B^{[n,m]} = (U^{[n,m]})^* U^{[n,m]} \left ( (U^{[n,m]})^* C^{[m]} U^{[n,m]} \right )^{-1} (U^{[n,m]})^* U^{[n,m]},
}
and $C^{[m]}$ is the Gram matrix for $\{ \psi_j \}^{m}_{j=1}$.  In particular, if $\{ \phi_j \}^{n}_{j=1}$ is an orthonormal basis for $\rT_n$,
\bes{
D_{n,m} = \frac{1}{\sigma_{\min}(U^{[n,m]})},\quad \sec ( \theta_{n,m} ) = \frac{1}{\sqrt{\lambda_{\min} ( B^{[n,m]} ) }},
}
where $\sigma_{\min}(U^{[n,m]})$ and $\lambda_{\min} ( B^{[n,m]} ) $ denote the minimal singular value and eigenvalue of the matrices $U^{[n,m]}$ and $B^{[n,m]}$ respectively.
}
\prf{
The proof of this lemma is similar to that of \cite[Lem.\ 2.13]{BAACHAccRecov}, and hence is omitted.
}

Although this lemma allows for computation of the reconstruction constant $C(F_{n,m})$ (recall that $C(F_{n,m}) =  \max \{ \sec (\theta_{n,m}) , D_{n,m} \}$ as a result of Corollary \ref{c:GSreconconst}), and therefore $\Theta(n;\theta)$ and $\Psi(m;\theta)$, it is somewhat inconvenient to have to compute both $D_{n,m}$ and $\sec (\theta_{n,m})$.  The latter, in particular, can be computationally intensive since it involves both forming and inverting the matrix $(U^{[n,m]})^* C^{[m]} U^{[n,m]}$.  However, recalling the bound $C(F_{n,m} ) \leq \max \{ 1 , \sqrt{c_2} \} D_{n,m}$, we see that stability and quasi-optimality can be ensured, up to the magnitude of $c_2$, by controlling the behaviour of $D_{n,m}$ only.  This motivates the computationally more convenient alternative
\bes{
\tilde{\Theta}(n ; \theta) = \min \left \{ m \in \bbN : D_{n,m} \leq \theta \right \},\quad n \in \bbN,\ \theta \in \left ( \frac{1}{\sqrt{c_1} \cos ( \theta_{\rT \rS} ) } , \infty \right ),
}
and likewise $\tilde{\Psi}(m;\theta)$.  Note that setting $m \geq \tilde \Theta (n;\theta)$ ensures a condition number of at worst $\theta$ and a quasi-optimality constant of at most $\max \{ 1 , \sqrt{c_2} \} \theta$.

\section{Optimality of generalized sampling}\label{s:optimality}
In the previous sections we provided an analysis of generalized sampling, which improved on \cite{BAACHShannon,BAACHAccRecov} by providing sharp bounds and establishing the connection between generalized sampling and certain oblique projections.  The purpose of this section is to address the question of optimality of generalized sampling; a topic which was not considered in either previous paper.  

We consider the following problem:

\prob{
\label{p:reconstruction}
Given the $m$ measurements $\{ \ip{f}{\psi_j}\}^{m}_{j=1}$ of an element $f \in \rH$, compute a reconstruction $\tilde{f}$ of $f$ from the subspace $\rT_n$.
}

Generalized sampling provides a (perhaps the most) straightforward solution to this problem -- namely, performing a least-squares fit of the data -- with stability and quasi-optimality being determined by the quantities $D_{n,m}$ and $\sec (\theta_{n,m})$.  An obvious question to pose is the following: can a different method outperform generalized sampling?  Our first answer to this question is given in the next section.

\subsection{Optimality amongst perfect methods}
\thm{
\label{t:optimality}
Suppose that $m,n \in \bbN$ are such that $D_{n,m} \neq 0$, where $D_{n,m}$ is given by \R{Dnm_def}.  Let $G_{n,m}$ be a method taking measurements $\{ \ip{f}{\psi_j}\}^{m}_{j=1}$ and giving a reconstruction $G_{n,m}(f) \in \rT_n$.  Suppose that $G_{n,m}$ is perfect in the sense of Definition \ref{d:perfect}.  Then, if the condition number $\kappa(G_{n,m})$ is defined in \R{condnumb}, we have 
\bes{
\kappa(G_{n,m}) \geq D_{n,m}.
}
In particular, if $F_{n,m}$ is the generalized sampling reconstruction, then $\kappa(G_{n,m}) \geq \kappa(F_{n,m})$.
}
\prf{
Since $G_{n,m}$ is perfect, we have $G_{n,m}(0) = 0$.  Setting $f = 0$ in \R{condnumb}, we notice that
\bes{
\kappa(G_{n,m}) \geq \lim_{\epsilon \rightarrow 0^+} \sup_{\substack{g \in \rH \\ \hat{g}^{[m]}  \neq 0}} \frac{\| G_{n,m}
(\epsilon g) \|}{\| \epsilon \hat{g}^{[m]} \|_{\ell^2}}.
}
Since $D_{n,m} \neq 0$, we have that $\hat{g}^{[m]}  \neq 0$ for $g \in \rT_n$ if and only if $g \neq 0$.  Thus, using the perfectness of $G_{n,m}$ once more,
\bes{
\kappa(G_{n,m}) \geq \lim_{\epsilon \rightarrow 0^+} \sup_{\substack{g \in \rT_n \\ g \neq 0}} \frac{\| G_{n,m}(\epsilon g) \|}{\| \epsilon \hat{g}^{[m]} \|} =  \lim_{\epsilon \rightarrow 0^+} \sup_{\substack{g \in \rT_n \\ g \neq 0}} \frac{\|g \|}{\|\hat{g}^{[m]} \|} = D_{n,m},
}
as required.  The second result follows from Corollary \ref{c:GSreconconst}.
}

This theorem, which is embarrassingly simple to prove, states the following: any perfect method for Problem \ref{p:reconstruction} must have a worse condition number than that of generalized sampling.  We remark that perfectness is not an unreasonable assumption in practice.  In particular, any method which is quasi-optimal (for fixed $n$ and $m$) is also perfect.  Indeed, quasi-optimality is equivalent to the condition
\be{
\label{cNM}
\| f - G_{n,m}(f) \| \leq \mu(G_{n,m}) \inf_{\phi \in \rT_n} \| f - \cQ_{n} f \|,\quad \forall f \in \rH,
}
for some $\mu(G_{n,m}) < \infty$ (the quasi-optimality constant).  The right-hand side vanishes for any $f \in \rT_n$, which implies perfectness of $G_{n,m}$.

One can also generalize Theorem \ref{t:optimality} somewhat to consider a larger class of methods.  Indeed, let $G_{n,m}$ be a method such that
\bes{
\| f - G_{n,m} f \| \leq \lambda \| f \|,\quad \forall f \in \rT_n,
}
for some $\lambda \in (0,1)$.  We refer to such methods as \textit{contractive}.  Note that perfect methods are a particular example of contractive methods with $\lambda = 0$.  Arguing as in the proof of Theorem \ref{t:optimality}, one can show that
\bes{
\kappa(G_{n,m}) \geq (1-\lambda) \kappa(F_{n,m}).
}
Hence, the condition number of generalized sampling can only possibly be improved by a factor of $(1-\lambda)$ when using a contractive method.

Beside the imposition of perfectness (or contractiveness), it may also appear at first sight that Theorem \ref{t:optimality} is also restrictive because it deals with a worst case scenario taken over the whole of $\rH$.  In practice, it may be the case that our interest does not lie with recovering all $f \in \rH$, but rather only those $f$ belonging to some subspace $\rU$ of $\rH$.  For example, $\rU$ could consist of functions with particular smoothness.  One may reasonably ask: is it possible to circumvent this bound if one restricts ones interest to only a small class of functions?  The answer is no.  If $\rT_n \subseteq \rU$, then one can just redefine the condition number $\kappa$ to be taken as a supremum over $\rU$, as opposed to $\rH$, and repeat the same argument.

It is also worth observing that Theorem \ref{t:optimality} can be weakened by considering the condition number $\kappa_{f}$ at a fixed $f \in \rH$ (we refer to \R{local_condition} for the definition of $\kappa_f$).  Indeed, it is clear from the proof that
\bes{
\kappa_f(G_{n,m}) \geq D_{n,m} \geq \kappa_{f}(F_{n,m}),\quad \forall f \in \rT_n.
}
In some applications, typically where one's interest lies with recovering only one fixed signal $f$, the local condition number is arguably more important.  Hence, the fact that condition numbers cannot be improved, even locally, demonstrates the importance of appropriately scaling $m$ with $n$.

One can also reformulate the conclusions of Theorem \ref{t:optimality} in terms of the stable sampling rate.  To this end, suppose that $G_{n,m}$ is any reconstruction method satisfying \R{cNM}, and let $\kappa(G_{n,m})$ and $\mu(G_{n,m})$ be its condition number and quasi-optimality constant respectively.  Define the reconstruction constant $C(G_{n,m}) = \max \left \{ \kappa(G_{n,m}) , \mu(G_{n,m}) \right \}$ in the standard way, and let
\bes{
\Theta_{G}(n;\theta) = \min \left \{ m \in \bbN : C(G_{n,m}) < \theta \right \},\quad n \in \bbN,
}
be the stable sampling rate for $G_{n,m}$.  If, given $n$, there does not exist an $m$ such that $C(G_{n,m}) < \theta$, then we set $\Theta_{G}(n;\theta) =\infty$.  In other words, for this $\theta$ and $n$, there is no number of samples $m$ which renders the reconstruction stable and quasi-optimal.  If this is not the case, then Theorem \ref{t:optimality} trivially gives that
\be{
\label{GSSR}
\Theta_{G}(n;\theta) \geq \Theta(n;\max \{ 1,\sqrt{c_2} \} \theta),
}
where $\Theta$ is the stable sampling rate for generalized sampling.  This result implies the following: up to a constant on the order of $\max \{ 1,\sqrt{c_2} \}$, any reconstruction requires at least the same number of samples as generalized sampling to guarantee a stable quasi-optimal reconstruction with constant $\theta$.  In applications (see \S \ref{s:examples}) one typically has that $\Theta(n;\theta) \sim c(\theta) g(n)$ for functions $c$ and $g$ with $c(\theta)$ decreasing as $\theta \rightarrow \infty$ and $g(n)$ increasing as $n \rightarrow \infty$.  In other words, the stable sampling rate $\Theta(n;\theta)$ grows asymptotically like $g(n)$ (typically $g(n) = n^{\alpha}$ for some $\alpha \geq 1$).  Hence, \R{GSSR} implies that no perfect method $G_{n,m}$ can have a stable sampling rate that grows at a slower rate than that of generalized sampling, although the constant can be slightly improved whenever the sampling frame has $c_2 >1$.

\subsection{An optimality result for problems with linear stable sampling rates}

Suppose that the stable sampling rate $\Theta(n;\theta)$ is linear in $n$  for a particular example of Problem \ref{p:reconstruction}.  This means that there is, up to a constant, a one-to-one correspondence between samples and reconstructed coefficients, which suggests that generalized sampling can only be outperformed by a constant factor in terms of the convergence of the reconstruction.  Another method for the problem (perfect or otherwise) might obtain a slightly smaller error, but the asymptotic \textit{rate} of convergence should be equal.  

This is formalized in the following theorem:
\thm{
\label{t:optimality_linear}
Let $\{ \psi_j \}_{j \in \bbN}$ be a frame for $\rH$, and let $\{ \rT_n \}_{n \in \bbN}$ a sequence of finite-dimensional subspaces satisfying \R{Tcond1} and \R{Tcond2}.  Suppose that the corresponding stable sampling rate $\Theta(n;\theta)$ is linear in $n$.  Let $f \in \rH$ be fixed, and suppose that there exists a sequence of mappings
\bes{
G_m : \{ \hat{f}_j \}^{m}_{j=1} \mapsto G_m(f) \in \rT_{\Psi_f(m)},
}
where $\Psi_f : \bbN \rightarrow \bbN$ with $\Psi_f(m) \leq c m$ for some $c > 0$.  Suppose also that there exist constants $c_1(f),c_2(f),\alpha_f > 0$ such that
\be{
\label{algconv}
c_1(f) n^{-\alpha_f} \leq \| f - \cQ_n f \| \leq c_2(f) n^{-\alpha_f},\quad \forall n \in \bbN.
}
Then, given $\theta \in \left (\frac{\max \{ 1 , \sqrt{c_2} \}}{\sqrt{c_1} \cos (\theta_{\rT \rS} ) },\infty \right )$, there exists a constant $c_f(\theta) > 0$ such that
\be{
\label{decay}
\| f - F_{\Psi(m;\theta),m}(f) \| \leq c_f(\theta) \| f - G_m(f) \|,\quad \forall m \in \bbN,
}
where $F_{n,m}$ corresponds to generalized sampling and $\Psi(m;\theta)$ is the stable reconstruction rate \R{SRR}.
}
\prf{
Since generalized sampling is quasi-optimal, 
\bes{
\| f - F_{\Psi(m;\theta),m}(f) \| \leq \theta \| f - \cQ_{\Psi(m;\theta)} f \|.
}
Using \R{algconv} we deduce that
\bes{
\| f - F_{\Psi(m;\theta),m}(f) \| \leq \theta \frac{c_2(f)}{c_1(f)} \left ( \frac{\Psi_f(m)}{\Psi(m;\theta)} \right )^{\alpha_f} \| f - \cQ_{\Psi_f(m)} f \|.
}
The orthogonal projection $\cQ_n f$ is the best approximation to $f$ from the subspace $\rT_n$.  Therefore
\bes{
\| f - F_{\Psi(m;\theta),m}(f) \| \leq \theta \frac{c_2(f)}{c_1(f)} \left ( \frac{\Psi_f(m)}{\Psi(m;\theta)} \right )^{\alpha_f} \| f - G_m(f) \|.
}
The result now follows from the fact that $\Psi(m;\theta) = \ord{m}$ and $\Psi_f(m) \leq c m$.
}

This theorem states that, in the case of a linear stable sampling rate, and for functions with algebraic decay of $\| f - \cQ_n f\|$, generalized sampling can only be improved upon by a constant factor.  As shown by \R{decay}, the error of generalized sampling decays at the same (or better) asymptotic rate as any other reconstruction method $G_m$.  Note that the stipulation of algebraic convergence \R{algconv} is reasonable in practice.  In the next section we shall see several examples for which this condition holds.

Unlike the case of generalized sampling, the method $G_m$ in the above theorem can depend in a completely nontrivial manner on the function $f$.  However, even with this added flexibility, this theorem shows that it is only possible to improve on generalized sampling by a constant factor.  An example of such a method is an oracle.  Suppose there was some method that, for a particular $f$ satisfying (\ref{algconv}), could recover the orthogonal projection $\cQ_m f$ \textit{exactly} (i.e.\ with no error) from $m$ samples.  The conclusion of the above corollary is that generalized sampling commits an error that is at worst a constant factor larger than that of this method.

\rem{
\label{r:nonlinear}
The fact that the stable sampling rate is linear is key to Theorem \ref{t:optimality_linear}.  In situations where $\Theta(n;\theta)$ is superlinear (for an example, see the next section), it is possible to devise methods, albeit unstable methods, with asymptotically faster rates of convergence. 
}

\section{Uniform resampling with generalized sampling}\label{s:examples}
\S \ref{s:introduction}--\ref{s:optimality} of this paper have considered generalized sampling in its abstract form.  We now consider its application to a particular problem, the so-called \textit{uniform resampling} (URS) problem.  As we show, generalized sampling leads to an improvement over the standard approach to this problem, which results in an ill-posed discrete reconstruction.

We first describe the URS problem in further detail.

\subsection{The uniform resampling problem}
In applications such as MRI, radio-astronomy and diffraction tomography \cite{GelbNonuniformFourier,RosenfeldURS1}, the URS problem addresses the question of how to recover the Fourier coefficients of a function $f \in \rL^2(-1,1)^d$ from nonuniformly spaced pointwise samples of its Fourier transform
\bes{
\hat{f}(\omega) = \frac{1}{2^{\frac{d}{2}}} \int_{(-1,1)^d} f(x) \E^{-\I \omega \pi x} \D x,\quad \omega \in \bbR.
}
This problem is important since typical sampling devices (such as MR scanners) are not best suited to acquire Fourier data in a uniform pattern (i.e.\ Fourier coefficients).  Indeed, it is often more convenient to acquire samples along interlacing spirals or radial lines, for example (see \cite{VanDeWalleEtAllURS} and references therein).  In uniform resampling, one seeks to compute Fourier coefficients from these nonharmonic Fourier samples, and then recover the image via a standard DFT.

Consider the case $d=1$, and let $\omega_{-n} < \omega_{-n+1} < \ldots < \omega_n$ be a set of $2n+1$ nonequispaced points at which $\hat{f}(\omega)$ is sampled.  The derivation of the standard URS reconstruction follows from the Shannon Sampling theorem \cite{RosenfeldURS1,RosenfeldURS2,GelbNonuniformFourier}.  Using this theorem, we have
\bes{
\hat{f}(\omega) = \sum_{k \in \bbZ} \hat{f}(k) \mathrm{sinc}(\omega-k),\quad \omega \in \bbR,
}
where the right-hand side converges uniformly, and therefore
\be{
\label{nonunifShannon}
\hat{f}(\omega_j) = \sum_{k \in \bbZ} \hat{f}(k) \mathrm{sinc}(\omega_j-k),\quad |j| \leq n.
}
Let $\alpha_k$, $k=-n,\ldots,n$ be the values $\alpha_k \approx \hat{f}(k)$ that we seek to compute from the samples $\{ \hat{f}(\omega_j) \}_{|j| \leq n}$.  It is natural to truncate \R{nonunifShannon} at level $n$, leading to
\bes{
\hat{f}(\omega_j) \approx \sum_{|k| \leq n} \alpha_k \mathrm{sinc}(\omega_j - k),\quad |j| \leq n.
}
Let $U^{[n,n]} \in \bbC^{2n+1,2n+1}$ be the matrix with $(j,k)^{\rth}$ entry $\mathrm{sinc}(\omega_j-k)$.  The URS method determines the vector $\alpha^{[n,n]} = \{ \alpha_k \}_{|k| \leq n}$ as the solution to the linear system
\be{
\label{URSmethod0}
U^{[n,n]} \alpha^{[n,n]} = \hat{f}^{[n]},
}
where $\hat{f}^{[n]} = \{ \hat{f}(\omega_j) \}_{|j| \leq n}$.

Suppose now that the finite collection $\{ \omega_j \}_{|j| \leq n}$ extends to an infinite set $\{ \omega_j \}_{j \in \bbZ}$ such that the system
\bes{
\psi_j(x) = \frac{1}{\sqrt{2}} \E^{\I \omega_j \pi x},\quad j \in \bbZ,
}
is a frame for $\rL^2(-1,1)$.  Let $\phi_j(x) = \frac{1}{\sqrt{2}} \E^{\I j \pi x}$, so that
\be{
\label{Fourierbasis}
\rT_n = \spn \left \{ \phi_j : |j| \leq n \right \},
}
is the space of trigonometric polynomials of degree $n$.  Then the URS method \R{URSmethod0} is nothing more than a specific instance of the consistent sampling framework described in \S \ref{s:consistent}.

It has been widely reported that the URS method \R{URSmethod0} may be very ill-conditioned in practice \cite{RosenfeldURS1,GelbNonuniformFourier}.  Various strategies have been applied to the linear system \R{URSmethod} to try to overcome this issue, with the most common involving first manually computing a singular value decomposition and then applying standard regularization techniques from the literature on discrete ill-posed problems \cite{RosenfeldURS1,GelbNonuniformFourier}.  However, this approach is both computationally expensive and sensitive to noise (see \cite{RosenfeldURS2} and references therein).  As a consequence, even in the presence of low noise, the resulting image can often be highly contaminated (much as in the examples presented in \S \ref{ss:consistfail}).

Although the effect of noise can be somewhat mitigated \cite{RosenfeldURS2}, it can never truly be removed, since the underlying discrete problem is ill-posed.  However, the interpretation of URS as an instance of consistent sampling means that this ill-posedness -- which is merely another instance of that seen in the consistent reconstructions of \S \ref{s:consistent} -- is completely artificial.    The key theorems presented in \S \ref{s:GS} demonstrate that by replacing \R{URSmethod0} with an overdetermined the least squares (i.e.\ generalized sampling)
\be{
\label{URSmethod}
U^{[n,m]} \alpha^{[n,n]} =\hat{f}^{[m]},
}
and by increasing $m$ suitably we will be able to obtain a numerically stable reconstruction.  Hence, rather than performing an intensive regularization on a discrete ill-posed problem, we discretize differently so as to obtain a well-posed discrete problem (recall the operator-theoretic interpretation of \S \ref{ss:optheory})

\rem{
Of course, Theorem \ref{t:optimality} states that if the uniform resampling \R{URSmethod0} is ill-conditioned (for a given $n$) then so will any other perfect method.  In other words, there is in essence no stable way to obtain $n$ Fourier coefficients from $n$ nonuniformly spaced Fourier samples.  Hence, increasing $m$ is not just a good way to proceed, in this sense it is the only possible way to obtain a stable reconstruction.

We remark also that overdetermined least squares of the form \R{URSmethod} has been used in the past for the uniform resampling problem.  However, it is still reported as resulting in an ill-conditioned problem \cite{RosenfeldURS1,GelbNonuniformFourier}.  This is unsurprising in the results of this paper:  $m$ needs to not only be larger than $n$ to ensure stability, but also above the critical threshold of the stable sampling rate $\Theta(n;\theta)$.
}

\rem{
There are a number of alternatives to uniform resampling, such as convolutional gridding techniques \cite{JacksonEtAlGridding,MeyerEtAllGridding,RosenfeldURS1,GelbNonuniformFourier}, which is quite popular in practice.  However, URS provides an optimal solution to the problem, and consequently often provides better results \cite{RosenfeldURS1} (in particular, it can lead to a significant decrease in artifacts \cite{VanDeWalleEtAllURS} over gridding).  Convolutional gridding has the advantage of being more efficient \cite{RosenfeldURS1} than the standard URS algorithm.  However, modifications such at the \textit{block uniform resampling} (BURS) \cite{RosenfeldURS1} possess comparable efficiency.
}

\subsection{Generalized sampling for the URS problem}
Provided one selects the parameter $m$ using the stable sampling rate for this problem, the key theorems of \S \ref{s:GS} demonstrate that \R{URSmethod} will be perfectly stable, as well as quasi-optimal.  It is therefore critical to determine $\Theta(n;\theta)$ in this instance.  Our main result below demonstrates that $\Theta(n;\theta)$ is linear in $n$ for (almost) all nonuniform sampling patterns arising as Fourier frames.  

First we require the following definition \cite{GelbHines2011Frames}:

\defn{
A sequence $\{ \omega_j \}_{j \in \bbZ}$ is a balanced sampling sequence if the following conditions hold:
\enum{
\item[(i)] $\Omega = \{ \E^{\I \omega_j \pi \cdot } : j \in \bbZ \}$ is a frame for $\rL^2(-1,1)$,
\item[(ii)] $\{ \omega_j \}_{j \in \bbZ}$ is $\delta$-separated, i.e.\ there exists $\delta>0$ such that $| \omega_j - \omega_k | > \delta$, $\forall j \neq k$,
\item[(iii)] $\{ \omega_j \}_{j \in \bbZ}$ is increasing, i.e.\ $\omega_{j} \leq \omega_{j+1}$, $\forall j \in \bbZ$,
\item[(iv)] $\{ \omega_j \}_{j \in \bbZ}$ is balanced, i.e.\ $\omega_j \geq 0$ if $j \geq 0$ and $\omega_j < 0$ if $j < 0$.
}
}
Note conditions (iii) and (iv) can always be guaranteed by reordering.  Condition (iv) is also reasonable in practice since sampling strategies are typically symmetric.  Although (ii) does not hold for all Fourier frames, we shall assume it for simplicity in the presentation that follows.  It is possible in what follows to derive a fully general result on the stable sampling rate for arbitrary Fourier frames using \cite[Lem.\ 2]{Jaffard} (see also \cite[Thm.\ 3]{GelbHines2011Frames}).  However, for simplicity we shall not do this.

\thm{
\label{t:URSSSR}
Suppose that $\{ \omega_j \}_{j \in \bbZ}$ is a balanced sampling sequence.  Then the stable sampling rate $\Theta(n;\theta) = \ord{n}$.  Specifically, let $\tau : \bbN \rightarrow (0,\infty)$ be given by $\tau(m) = \min \{ \omega_m , - \omega_{-m} \}$, and define $\tau^{-1} : (0,\infty) \rightarrow \bbN$ by
\bes{
\tau^{-1}(c) = \min \{ m \in \bbN : \tau(m) > c \}.
}
Then $\tau^{-1}(c) < \lceil \frac{c}{\delta} \rceil $, $\forall c > 0$, and we have the upper bound
\bes{
\Theta(n;\theta) \leq \tau^{-1}\left ( \frac{g(\theta)}{g(\theta)-1} + \frac{g(\theta)+1}{g(\theta)-1} n \right ),
}
where $g(\theta) = \exp \left ( \pi^2 \delta ( c_1 - \max \{ 1 , c_2 \} \theta^{-2} ) \right )$.
}
\prf{
Let $\cP_m g = \sum_{|j| \leq m} \ip{g}{\psi_j} \psi_j$, and suppose that $\phi \in \rT_n$ is arbitrary.  Then
\be{
\label{decay2}
\ip{\cP_m \phi}{\phi} = \ip{\cP \phi}{\phi} - \ip{(\cP - \cP_m) \phi}{\phi} \geq c_1 \| \phi \|^2 - \ip{(\cP - \cP_m) \phi}{\phi}.
}
Let $\phi = \sum_{|j| \leq n} \alpha_j \phi_j$ so that $\| \phi \| = \| \alpha \|_{\ell^2}$.  Since
\bes{
\ip{(\cP - \cP_m) f}{g} \leq \sqrt{\ip{(\cP - \cP_m) f}{f}} \sqrt{\ip{(\cP - \cP_m) g}{g}},\quad \forall f,g \in \rL^2(-1,1),
}
it follows that
\ea{
\ip{(\cP - \cP_m) \phi}{\phi} &= \sum_{|j|,|k| \leq n} \alpha_j \overline{\alpha_k} \ip{(\cP-\cP_m) \phi_j}{\phi_k} \nn
\\
& \leq \left ( \sum_{|j| \leq n} |\alpha_j| \sqrt{ \ip{(\cP-\cP_m) \phi_j}{\phi_j} } \right )^2 \nn
\\
& \leq \| \phi \|^2 \sum_{|j| \leq n} \ip{(\cP-\cP_m) \phi_j}{\phi_j} . \label{phi_m_err}
}
Let us suppose that $m$ is sufficiently large so that $| \omega_j | > n$ for $|j | > m$ (i.e.\ $m > \frac{n}{\delta}$).  Note that
\bes{
|\ip{\phi_k}{\psi_j}| \leq \frac{1}{\pi | \omega_j - k |},\quad |j| > m, \quad |k|\leq n.
}
Therefore
\bes{
\ip{(\cP - \cP_m)\phi_k}{\phi_k} \leq \sum_{|j|>m} \frac{1}{\pi^2 | \omega_j-k|^2} = \frac{1}{\pi^2} \sum_{j > m} \frac{1}{| \omega_j - k |^2} + \frac{1}{\pi^2} \sum_{j>m } \frac{1}{| \omega_{-j} - k |^2} .
}
Consider the first sum.  We have
\bes{
 \sum_{j > m} \frac{1}{| \omega_j - k |^2}  \leq \frac{1}{\delta} \int^{\infty}_{\omega_m} \frac{1}{(\omega - k )^2} \D \omega =  \frac{1}{\delta(\omega_m - k) }.
}
Using a similar estimate for the other sum, we obtain
\bes{
\ip{(\cP - \cP_m)\phi_k}{\phi_k} \leq \frac{1}{\pi^2 \delta} \left (  \frac{1}{\omega_m - k } +  \frac{1}{k-\omega_{-m} } \right ).
}
Substituting this into \R{phi_m_err}, we obtain
\bes{
\ip{(\cP - \cP_m) \phi}{\phi}  \leq \frac{ \| \phi \|^2}{\pi^2 \delta} \sum_{|k| \leq n} \left (  \frac{1}{\omega_m - k } +  \frac{1}{k-\omega_{-m} } \right ).
}
Notice that 
\bes{
\sum_{|k| \leq n}  \frac{1}{\omega_m - k } \leq \int^{n+1}_{-n} \frac{1}{\omega_m - x} \D x = \ln \left ( \frac{\omega_m + n}{\omega_m - n- 1} \right ).
}
Likewise
\bes{
\sum_{|k| \leq n}  \frac{1}{k-\omega_{-m}} \leq \int^{n}_{-n-1} \frac{1}{x-\omega_{-m}} \D x = \ln \left ( \frac{n-\omega_{-m}}{-\omega_{-m}-n-1} \right ).
}
Hence
\bes{
 \ip{(\cP - \cP_m) \phi}{\phi}  \leq \frac{\| \phi \|^2}{\pi^2 \delta} \ln \left ( \frac{\tau(m) + n}{\tau(m) - n- 1} \right ).
}
Combining this with \R{decay2}, we obtain
\bes{
D^{-2}_{n,m} =\inf_{\substack{\phi \in \rT_n \\ \| \phi \|=1}} \ip{\cP_m \phi}{\phi} \geq c_1 - \frac{1}{\pi^2 \delta} \ln \left ( \frac{\tau(m) + n}{\tau(m) - n- 1} \right ).
}
Recall that $C(F_{n,m}) \leq \max \{ 1, \sqrt{c_2} \} D_{n,m}$.  Hence $C(F_{n,m}) < \theta$ provided 
\bes{
c_1 - \frac{1}{\pi^2 \delta} \ln \left ( \frac{\tau(m) + n}{\tau(m) - n- 1} \right ) > \max \{ 1 , c_2 \} \theta^{-2}.
}
Rearranging, we obtain
\bes{
\frac{\tau(m) + n}{\tau(m) - n- 1} < g(\theta) \quad \Longleftrightarrow \quad \tau(m) > \frac{g(\theta)}{g(\theta)-1} + \frac{g(\theta)+1}{g(\theta)-1}  n,
}
and this gives the result (note that this condition implies that $| \omega_j | > n$, $|j| > m$, which was the assumption made for the above analysis).  Note also that $\tau(m) = \min \{ \omega_{m} , - \omega_{-m} \} \geq \delta m$.  Thus $\tau^{-1}(c) \leq \lceil \frac{c}{\delta} \rceil $.  This completes the proof.
}

\subsection{Numerical results}
We now give numerical results for generalized sampling for applied to this problem.  We consider the following three sequences
\bes{
(a): \omega_{j} = \frac{1}{2} j,\quad (b): \omega_{j} = \frac{1}{4} j,\quad (c): \omega_{j} = \frac{1}{4} j + \nu_{j},
}
where in the last case $\nu_{j} \in (-\frac15,\frac15)$ is chosen uniformly at random. Note that all three sequences are frames for $\rL^2(-1,1)$ \cite{christensen2003introduction}.

\begin{figure}
\begin{center}
$\begin{array}{ccc}
 \includegraphics[width=4.75cm]{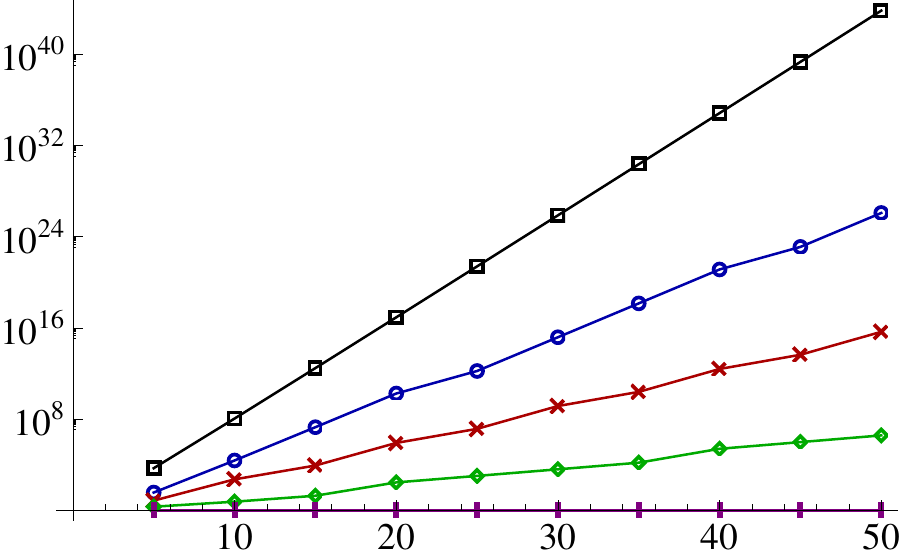}   &    \includegraphics[width=4.75cm]{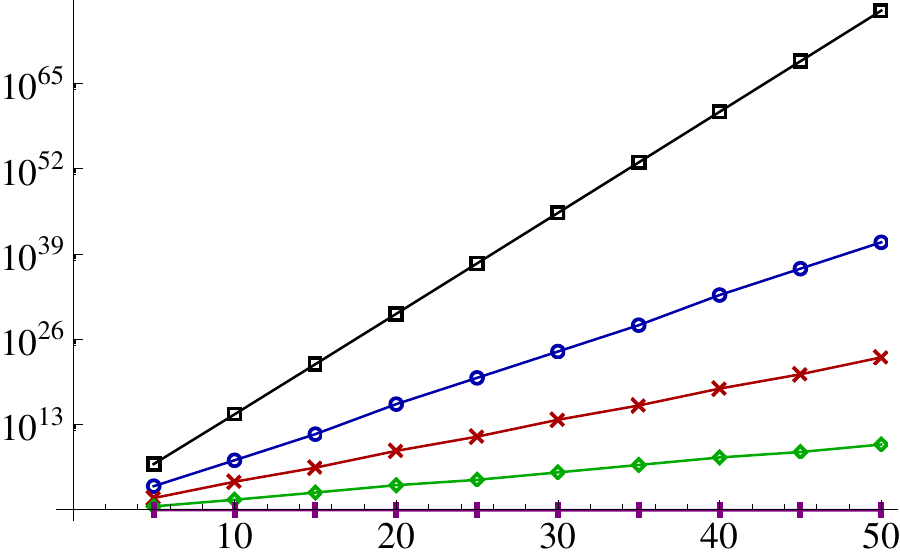} &    \includegraphics[width=4.75cm]{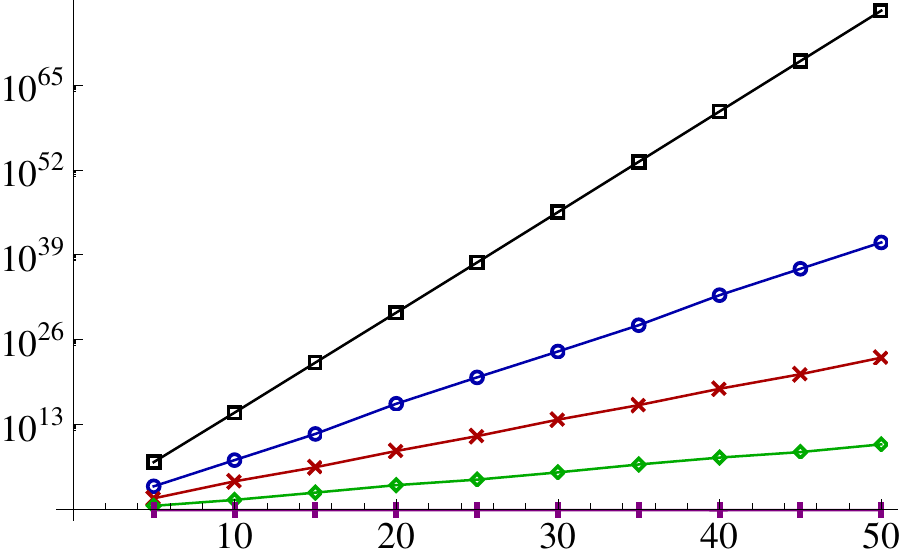} \\
 (a) & (b) & (c)
\end{array}$
\caption{\small The quantity $D_{n,c n}$ against $n$ for the generalized sampling applied to the uniform resampling problem for the frames (a)--(c).  The values $c=1,\frac{5}{4},\frac{3}{2},\frac{7}{4},2$ were used for (a) and $c=1,\frac{7}{4},\frac{5}{2},\frac{13}{4},4$ for (b) and (c).}  \label{f:URSReconConst}
\end{center}
\end{figure}

 In Figure \ref{f:URSReconConst} we plot the quantity $D_{n,m}$ with various linear scalings of $m$ with $n$.  As is evident, this constant is exponentially large when $m=n$ (i.e.\ consistent sampling), and it remains exponentially large when $m = c n$ for small $c$ below a certain threshold.  However, as $c$ increases the rate of exponential growth decreases, and once $c$ is sufficiently large, there is no exponential growth at all. 

\begin{figure}
\begin{center}
$\begin{array}{ccc}
 \includegraphics[width=4.75cm]{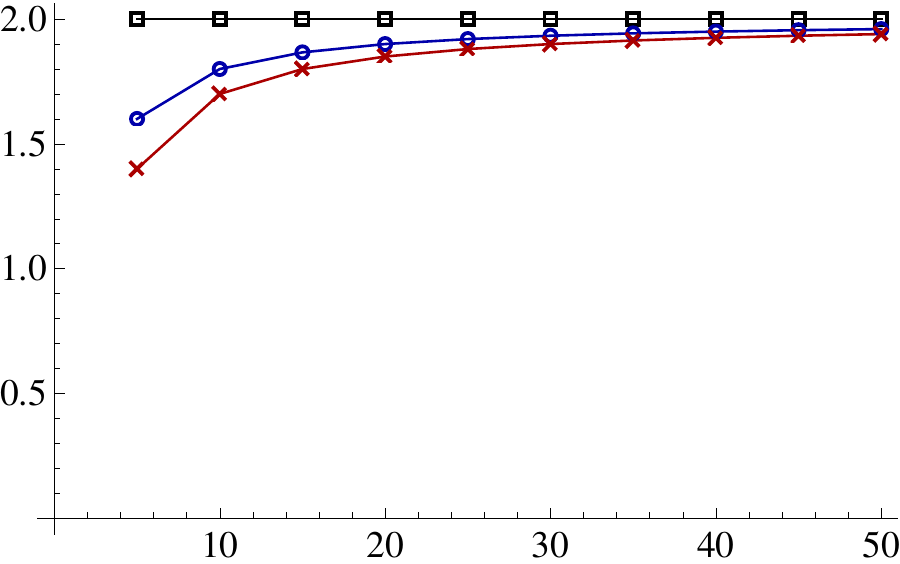}   &    \includegraphics[width=4.75cm]{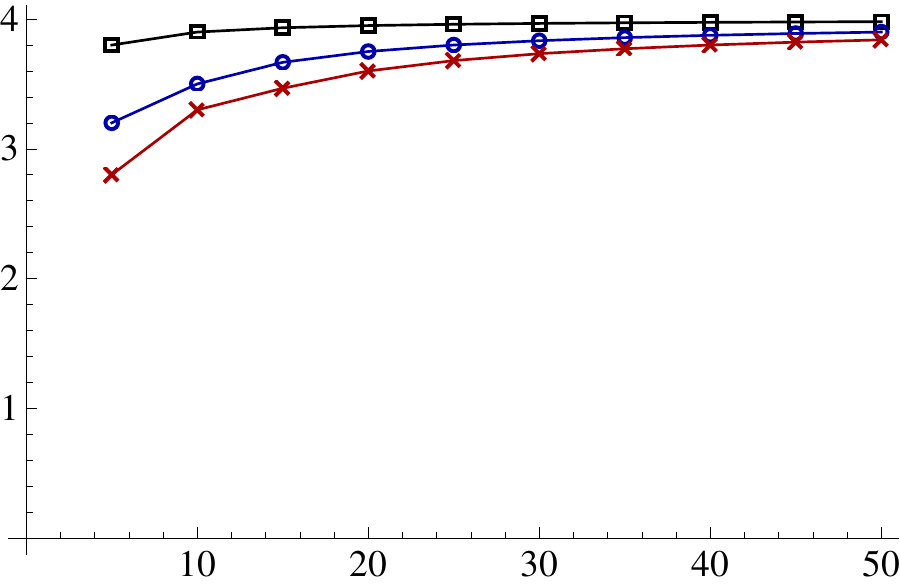} &    \includegraphics[width=4.75cm]{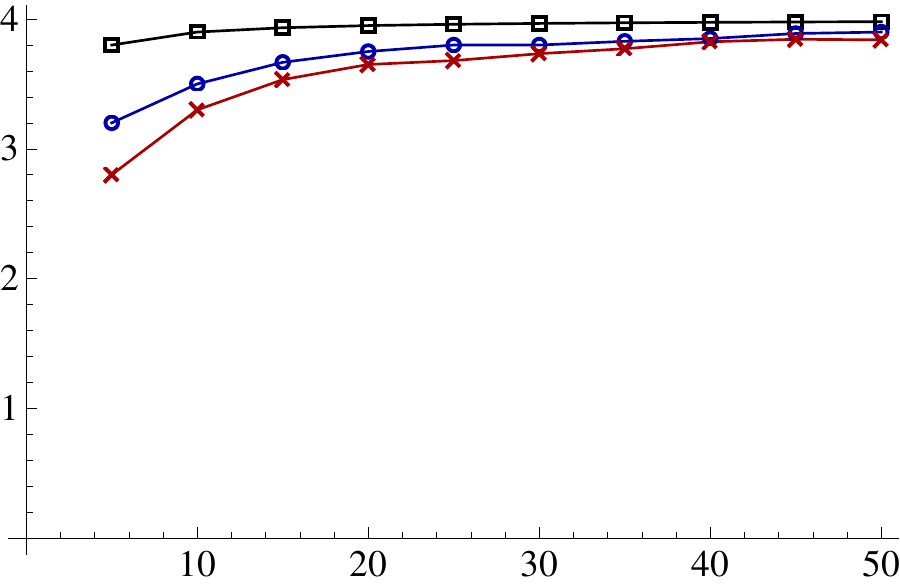} \\
 (a) & (b) & (c)
\end{array}$
\caption{\small The stable sampling rate $\tilde \Theta(n;\theta)$, scaled by $n^{-1}$, for the frames (a)--(c), where $\theta = \frac{5}{4},10,50$.}  \label{f:URSSSR}
\end{center}
\end{figure}

To determine the critical $c$ for which the reconstruction constant is bounded, we compute the stable sampling rate.  This is shown in Figure \ref{f:URSSSR}.  For the frame (a) this critical value is roughly $2$, whereas for (b) and (c) it is approximately $4$.  Moreover, the closeness of the graphs indicates that one only needs to exceed this critical value by a very small amount to get an extremely good reconstruction constant.

\begin{figure}
\begin{center}
$\begin{array}{ccc}
 \includegraphics[width=4.75cm]{Diagrams/Ex1Noise} & \hspace{2pc} & \includegraphics[width=4.75cm]{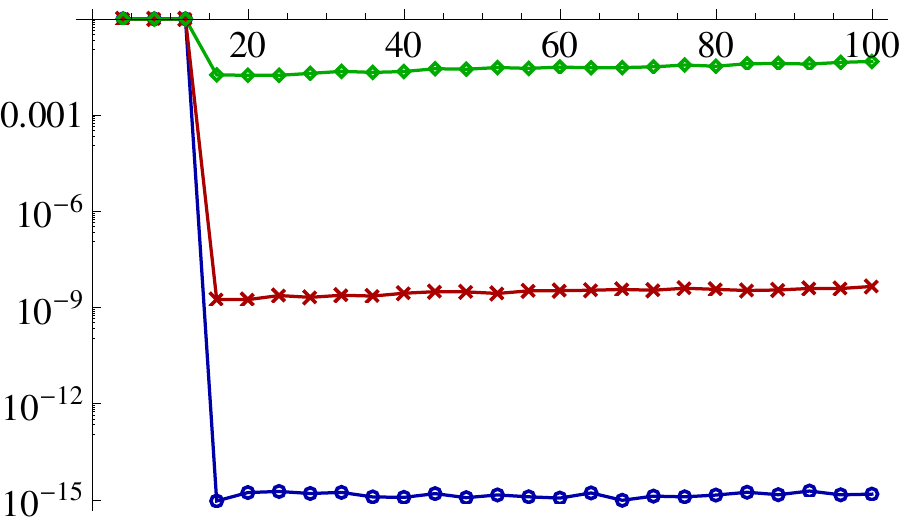}   
\end{array}$
\caption{\small The error $\| f - \tilde{f}_{n,m} \|$ against $m$ where $n=m$ (left) and $n=\frac12 m$ (right), $f(x) = \frac{1}{\sqrt{2}} \E^{8 \I \pi x}$, and $\tilde{f}_{n,m}$ is computed from noisy data $\{ \hat{f}_{j} + \eta_j \}_{|j| \leq n}$, where $| \eta_j | \leq \eta$ is chosen uniformly at random with $\eta = 0,10^{-9},10^{-2}$ (circles, crosses and diamonds respectively).  The sampling frame (a) was used.}  \label{f:URSNoise}
\end{center}
\end{figure}

To illustrate the effectiveness of generalized sampling for this problem, in Figure \ref{f:URSNoise} we consider the reconstruction from noisy data.  As is evident, when $m=n$, noise is amplified by around $10^{15}$.  However, double oversampling, as suggested in Figure \ref{f:URSSSR}, renders the reconstruction completely stable: the overall reconstruction error is on the order of the magnitude of the noise.

\subsection{Alternatives to uniform resampling}
The goal of uniform resampling is to recover the Fourier coefficients of the unknown function $f$ from its nonuniform Fourier samples.  However, it is well known that images and signals (which are typically nonperiodic) are poorly represented by their Fourier series.  Although the Fourier series converges (due to the Shannon Sampling Theorem), the rate is often intolerably slow and the finite series is polluted by Gibbs oscillations.

However, there is no reason besides familiarity to actually compute Fourier coefficients from nonuniform Fourier samples.  With generalized sampling one is able to reconstruct in \textit{any} subspace $\rT_n$; in particular, one which is better suited to the particular function.  Thus, provided such a subspace is known, we are able to obtain a better reconstruction over the classical Fourier series. 

In this final section we consider briefly two alternative choices for $\rT_n$ besides the Fourier space \R{Fourierbasis}.  The first is a spline space of piecewise polynomial functions of fixed degree $d \in \bbN$:
\be{
\label{splinespace}
\rT_n = \left \{ \phi \in \rC^{d-1}[-1,1] : \phi |_{[\frac{j}{n},\frac{j+1}{n})} \in \bbP^d,\quad j=-n,\ldots,n-1 \right \},\quad n \in \bbN.
}
Note that the sequence of orthogonal projections $\cQ_n f$ of a function $f \in \rC^d[-1,1]$ converge to $f$ at the rate $n^{-d-1}$, without the assumption of periodicity of $f$.  Conversely, the Fourier series (i.e.\ the URS reconstruction with $\rT_n$ given by \R{Fourierbasis}) converges like $n^{-\frac12}$ when $f$ is nonperiodic.  Hence, the spaces $\rT_n$ are better suited for moderately smooth and nonperiodic functions.

The second choice for $\rT_n$ is the polynomial space
\be{
\label{polyspace}
\rT_n = \bbP_n.
}
Observe that if $f$ is smooth, i.e.\ $f \in \rC^{\infty}[-1,1]$, then $\cQ_{n} f$ converges faster than any power of $n^{-1}$.  Hence, this space is particularly well suited for smooth functions.  Note that the use of this space for uniform Fourier samples $\omega_j = j$ was extensively discussed in \cite{BAACHAccRecov}.

\begin{table}
\centering 
\begin{tabular}{c|cccc}
$n$ &  8 & 16 &  32 & 64    \\ \hline 
(a)  & 8.33e12 & 7.24e25 & 2.17e52 & 3.40e105    \\
(b)  & 1.51e5 & 3.14e12 & 1.37e25 &  2.69e51  \\
(c) & 2.54e6 & 9.53e11 &  1.65e23  &  3.98e45 \\
\end{tabular} 
\quad
\begin{tabular}{c|cccc}
$n$ &  8 & 16 &  32 & 64    \\ \hline 
(a)  &  3.14e11 & 6.68e23 &  2.25e49 & 2.16e104   \\
(b) & 4.95e4 & 3.84e10 & 5.32e22 & 1.41e48   \\
(c) & 2.29e2 & 4.19e4 & 2.11e9   &  1.30e19 \\
\end{tabular} 
 \caption{\small The quantities $D_{n,n}$ (left) and  $D_{n,2n}$ (right) for the spaces (a): \R{splinespace} ($d=2$), (b): \R{splinespace} ($d=4$) and (c): \R{polyspace}.  The sampling frequencies are $\omega_j = \frac14 j + \nu_j$, where $\nu_j \in (-\frac15,\frac15)$ was chosen uniformly at random.}
  \label{ReconConstTable}
\end{table}

As one might expect, both the spaces \R{splinespace} and \R{polyspace} lead to instability in the corresponding consistent reconstruction.  This is shown in Table \ref{ReconConstTable}: in both cases, the constant $D_{n,n}$ is exponentially large in $n$.  Nonetheless, such instability can be overcome by sampling at the stable sampling rate.  Although we shall not do it in this paper (for the sake of brevity), it is possible to prove that the stable sampling rate is linear $\Theta(n;\theta) = \ord{n}$ for the spaces \R{splinespace} for any fixed $d \in \bbN$, and quadratic $\Theta(n;\theta) = \ord{n^2}$ for \R{polyspace}.  Note that a similar result for the latter in the case of uniform Fourier samples was shown previously in \cite{BAACHAccRecov} and \cite{hrycakIPRM}.

Instead, we now illustrate the advantage gained from exploiting these different reconstruction spaces.  In Tables \ref{Comparison1} and \ref{Comparison2} we give numerical results for the three different spaces considered.  In each case, the parameter $m$ (the number of samples) was fixed and $n$ chosen so that the quantity $D_{n,m} \leq 4$.  As can be seen in Table \ref{Comparison1}, the Fourier space \R{Fourierbasis} is particularly well suited for periodic functions, and outperforms both the spline \R{splinespace} and polynomial \R{polyspace} spaces.  However, the situation changes completely when the function to be reconstructed is not periodic.  In Table \ref{Comparison2} we see that the polynomial space \R{polyspace} gives the best reconstruction, followed by the spline space \R{splinespace}.  The URS reconstruction, which uses the Fourier space \R{Fourierbasis}, suffers from the Gibbs phenomenon and thus exhibits only low accuracy.

\begin{table}
\centering 
\begin{tabular}{c|cccc}
$m$ & Fourier & splines ($d=2$) & splines ($d=4$) & polynomials   \\ \hline 
32 & 2.80e-2 & 3.64e-3 & 5.95e-4 & 1.42e-2 \\
64 &  1.49e-4 &4.15e-4 & 3.05e-5 & 3.23e-3 \\
128 & 1.04e-11 & 4.89e-5 & 4.76e-7 & 3.01e-5\\
256 & 3.30e-15 & 6.04e-6 & 1.33e-8 & 5.43e-8\\
512 & 4.06e-15 &  7.55e-7 & 4.10e-10 &  5.02e-14 \\
1024 & 4.79e-15 & 9.45e-8 & 1.27e-11& 4.86e-14 \\
\end{tabular} 
 \caption{\small The error $\| f - \tilde{f}_{n,m} \|$ for the smooth and periodic function $f(x) = \sin 3 \pi x + 2 \E^{\frac{20}{\pi^2} ( \cos 2 \pi x - 4 \cos \pi x - 5)}$, where the reconstruction space $\rT_n$ is the Fourier space \R{Fourierbasis}, the spline space \R{splinespace} with $d=2,4$, or the polynomial space \R{polyspace}.  The sampling frequencies are given by $\omega_j = \frac14 j + \nu_j$, where $\nu_j \in (-\frac15,\frac15)$ was chosen uniformly at random.  The parameter $n$ was chosen so that $D_{n,m} \leq 4$.}
  \label{Comparison1}
\end{table}

\begin{table}
\centering 
\begin{tabular}{c|cccc}
$m$ & Fourier & splines ($d=2$) & splines ($d=4$) & polynomials   \\ \hline 
32 & 3.13e-1 & 3.12e-2 & 4.72e-3 & 3.62e-2  \\
64 &  1.53e-1 & 4.79e-3 & 1.80e-4 & 1.56e-3 \\
128 & 8.65e-2 & 5.16e-4 & 1.31e-5 &  1.79e-6  \\
256 & 9.27e-2 & 5.96e-5 & 4.84e-7 &  4.41e-11 \\
512 & 6.21e-2 & 7.14e-6 & 1.32e-8 & 4.33e-14 \\
1024 & 2.50e-2 & 8.82e-7 & 3.94e-10 & 4.19e-14
\end{tabular} 
 \caption{\small The error $\| f - \tilde{f}_{n,m} \|$ for the smooth function $f(x) = \sin 10 x + 2 \E^{20 (x^2 - 1)}$, where the reconstruction space $\rT_n$ is the Fourier space \R{Fourierbasis}, the spline space \R{splinespace} with $d=2,4$, or the polynomial space \R{polyspace}.  The sampling frequencies are given by $\omega_j = \frac14 j + \nu_j$, where $\nu_j \in (-\frac15,\frac15)$ was chosen uniformly at random.  The parameter $n$ was chosen so that $D_{n,m} \leq 4$.}
  \label{Comparison2}
\end{table}

\subsection{On optimality}\label{ss:splinepolyopt}

In view of the numerical results for the constant $D_{n,m}$ (Figure \ref{f:URSReconConst}), Theorem \ref{t:optimality} demonstrates that any perfect method for the uniform resampling problem with be exponentially unstable, unless the stable sampling rate is adhered to.  Moreover, since the stable sampling rate is linear in this case (Theorem \ref{t:URSSSR}), Theorem \ref{t:optimality_linear} also applies in this instance.  Hence, for periodic functions of finite smoothness (i.e.\ functions for which \R{algconv} holds), one cannot outperform generalized sampling by more than a constant factor regardless of the method.

The same conclusions also hold in the case of the spline spaces \R{splinespace}, in view of the numerics in Table \ref{ReconConstTable}.  For the polynomial space \R{polyspace}, however, Theorem \ref{t:optimality_linear} does not apply, since the stable sampling rate is quadratic.  Hence, it is in theory possible to outperform generalized sampling in terms of the asymptotic rate of convergence.  Nonetheless, it transpires that this cannot be done in this case without compromising stability.  For a more thorough analysis of stability and convergence for this reconstruction problem we refer the reader to \cite{AdcockHansenShadrinStabilityFourier}.

\section*{Acknowledgements}
The authors would like to thank Yonina Eldar, Hagai Kirschner, Nilima Nigam, Alexei Shadrin and Michael Unser for their helpful discussions and comments.

\bibliographystyle{abbrv}
\small
\bibliography{OptimalityRefs}

\end{document}